\documentclass{article}
\usepackage{amsmath,amsthm,amsfonts,eucal,amssymb}

\numberwithin{equation}{section}

\def\ca{{\mathcal A}}
\def\cb{{\mathcal B}}

\def\cd{{\mathcal D}}

\def\cf{{\mathcal F}}

\def\ch{{\mathcal H}}
\def\cai{{\mathcal I}}

\def\ck{{\mathcal K}}
\def\cl{{\mathcal L}}
\def\cam{{\mathcal M}}

\def\cp{{\mathcal P}}

\def\bn{{\mathbb N}}
\def\br{{\mathbb R}}

\def\a{\alpha}
\def\b{\beta}
\def\g{\gamma}        
\def\d{\delta}        \def\D{\Delta}
\def\eps{\varepsilon} 
     
\def\z{\zeta}
\def\th{\vartheta}

\def\l{\lambda}       \def\La{\Lambda}
\def\m{\mu}
\def\n{\nu}

\def\r{\rho}
\def\s{\sigma}
     \def\S{\Sigma}
\def\t{\tau}

\def\f{\varphi}
\def\o{\omega}        \def\O{\Omega}

\def\itm#1{\item[($#1$)]}

\DeclareMathOperator{\Lim}{Lim}

\DeclareMathOperator{\Tr}{Tr}
\DeclareMathOperator{\tr}{tr}
\DeclareMathOperator{\diam}{diam}
\DeclareMathOperator{\vol}{vol}
\DeclareMathOperator{\ran}{ran}

\def\ubd#1{\overline{d_{B}}(#1)} 
\def\lbd#1{\underline{d_{B}}(#1)}

\def\cpt{\ck(\ch)}

\def\subc{\underline{\d}}
\def\supc{\overline{\d}}
\def\subd{\underline{d}}
\def\supd{\overline{d}}
\def\subg{\underline{g}}
\def\supg{\overline{g}}

\def\sinc{S^{\uparrow}}
\def\sdec{S^{\downarrow}}
\def\triple{{\mathbb A}}

\newtheorem{Thm}{Theorem}[section]
\newtheorem{Cor}[Thm]{Corollary}
\newtheorem{Prop}[Thm]{Proposition}
\newtheorem{Lemma}[Thm]{Lemma}
\theoremstyle{definition}
\newtheorem{Dfn}[Thm]{Definition}

\theoremstyle{remark}
\newtheorem{rem}[Thm]{Remark} 
\newtheorem{ack}{Acknowledgement} 

\begin{document}
\title{\huge  Dimensions and singular traces   
 for spectral triples, \\
 with applications to fractals}
\author{Daniele Guido$^{*}$, Tommaso Isola
\thanks{Supported in part by GNAMPA and MIUR}
}
\date{November 12, 2002}
\markboth{Dimension, spectral triples and fractals}
{Dimension, spectral triples and fractals}
\maketitle
\bigskip\bigskip\noindent
Dipartimento di Matematica, Universit\`a di Roma ``Tor
Vergata'', I--00133 Roma, Italy. E-mail: {\tt guido@mat.uniroma2.it, 
isola@mat.uniroma2.it}

\begin{abstract}
	Given a spectral triple $(\ca,\ch,D)$, the functionals on $\ca$ of 
	the form $a\mapsto\t_{\o}(a|D|^{-\a})$ are studied, where 
	$\t_{\o}$ is a singular trace, and $\o$ is a generalised limit.  
	When $\t_{\o}$ is the Dixmier trace, the unique exponent $d$ giving rise 
	possibly to a non-trivial functional is called Hausdorff 
	dimension, and the corresponding functional the ($d$-dimensional) 
	Hausdorff functional.
	
	It is shown that the Hausdorff dimension $d$ coincides with the 
	abscissa of convergence of the zeta function of $|D|^{-1}$, and 
	that the set of $\a$'s for which there exists a singular trace 
	$\t_{\o}$ giving rise to a non-trivial functional is an interval 
	containing $d$.  Moreover, the endpoints of such traceability 
	interval have a dimensional interpretation.  The corresponding 
	functionals are called Hausdorff-Besicovitch functionals.
	
	These definitions are tested on fractals in $\br$, by computing 
	the mentioned quantities and showing in many cases their 
	correspondence with classical objects.  In particular, for 
	self-similar fractals the traceability interval consists only of 
	the Hausdorff dimension, and the corresponding 
	Hausdorff-Besicovitch functional gives rise to the Hausdorff 
	measure.  More generally, for any limit fractal, the described 
	functionals do not depend on the generalized limit $\o$.
\end{abstract}

 \newpage

\setcounter{section}{-1}

\section{Introduction.}\label{sec:zeroth}
The concept of spectral triple, introduced by Alain Connes as a 
framework for noncommutative geometry \cite{Co}, is wide enough to 
describe non smooth, or even fractal spaces.

While further axioms can be, and have been, added to describe the 
noncommutative analogues of Riemannian (or spin) manifolds (see 
\cite{GBVa} and references therein), spectral triples have been 
attached to fractals in $\br$ and to quasi-circles, and, using the 
Hausdorff dimension as the exponent for the infinitesimal length 
element, and the Dixmier trace, Connes, resp.  Connes-Sullivan proved 
that one obtains the Hausdorff measure for Cantor sets, resp.  
quasi-circles (cf.  \cite{Co}, IV.3).

Our aim in this paper is twofold. On the one hand we intend to show 
how spectral triples can provide a framework for noncommutative 
Hausdorff-Besicovitch theory. On the other hand we investigate how the 
noncommutative quantities we introduce give back classical known 
quantities, or even produce new ones, when applied to spectral 
triples associated to fractals.

Let us recall that a spectral triple $(\ca,\ch,D)$ consists of a 
$^{*}$-algebra $\ca$ acting on a Hilbert space $\ch$ and of a 
selfadjoint operator $D$ with compact resolvent, the Dirac operator, 
such that $[D,a]$ is bounded for any $a\in\ca$.

Concerning the first point, we define, as suggested by the examples of 
Connes, the $\a$-dimensional Hausdorff functional as the functional 
$a\mapsto\Tr_{\o}(a|D|^{-\a})$, where $\Tr_{\o}$ denotes the Dixmier 
trace, namely the singular trace summing logarithmic divergences.

Once Hausdorff functionals are defined, the Hausdorff dimension of a 
spectral triple is easily defined (cf.  Definition 
\ref{Dfn:dimensions}).  We show moreover that such dimension is equal 
to the abscissa of convergence of the zeta function of $|D|^{-1}$, 
$\zeta_{|D|^{-1}}(s)=\tr|D|^{-s}$.  This result turns out to be a 
useful formula for computing the Hausdorff dimension, and also opens 
the way to general Tauberian formulas for singular traces \cite{Taub}, 
see also \cite{CPS} for related results.

However our main goal here is to enlarge the class of geometric 
measures, in the same way as Besicovitch measures generalize Hausdorff 
measures.

Let us recall that the Hausdorff-Besicovitch measures replace the 
power law for the volume of the balls of the Hausdorff measures with a 
general infinitesimal behaviour.

Following the idea that the powerlike asymptotics for $\m_{n}(|D|)$, 
which give rise to non trivial logarithmic Dixmier traces, are the 
noncommutative counterpart of the powerlike asymptotics of the volume 
of the balls with small radius, which corresponds to some nontrivial 
Hausdorff measure, it is clear that, in order to define 
Hausdorff-Besicovitch functionals, we have to pass from the 
logarithmic singular trace to a general singular trace.

The trace theorem of Connes (\cite{Co}, IV.2) shows that the 
logarithmic Dixmier trace produces a trace functional on the 
C$^{*}$-algebra $\overline\ca$, corresponding to the Riemann volume on 
manifolds.  The proof of this theorem given in \cite{CiGS1} however, 
works for any positive trace functional whose domain contains the 
principal ideal generated by $|D|^{-1}$, which is singular, namely 
vanishes on finite rank operators, and for which the H\"older 
inequality holds.

We show here that H\"older inequality holds for any singular trace, 
possibly up to a constant (cf.  Appendix), namely any such trace 
produces a Hausdorff-Besicovitch functional on $\overline\ca$.

Then a new question arises: given a spectral triple $(\ca,\ch,D)$,
characterize the set $\{\a>0: |D|^{-\a}$ is singularly traceable $\}$,
namely describe when nontrivial Hausdorff-Besicovitch functionals can
be produced.  We give a complete answer to this question, namely prove
that such set is a relatively closed interval in $(0,\infty)$, whose
endpoints $\subc$ and $\supc$ coincide with the Matuszewska indices of
the eigenvalue function of $|D|^{-1}$, and satisfy $\subc\leq
d\leq\supc$, where $d$ denotes the Hausdorff dimension of the spectral
triple.

This means in particular that when $d$ is finite nonzero, it gives 
rise to a nontrivial Hausdorff-Besicovitch functional (see Theorem 
\ref{Thm:trace}).  Besides, when the traceability exponent is unique, 
it has to coincide with the Hausdorff dimension.  We remark that the 
singular trace associated with $d$ is not necessarily logarithmic 
Dixmier, indeed $|D|^{-d}$ may also be trace class.  The existence of 
a nontrivial Hausdorff-Besicovitch functional on $|D|^{-d}$ therefore 
depends essentially on the fact that all singular traces are allowed, 
not only those vanishing on $\cl^{1}$ (cf.  \cite{AGPS}).

The quantities $\subc$, $\supc$ mentioned above exhibit some 
dimensional behaviour (cf.  Theorem \ref{sumdim}), therefore we shall 
call them minimal, resp.  maximal, dimension of the spectral triple.  
In order to understand if these quantities, which have been introduced 
in a purely noncommutative fashion, have a commutative counterpart, we 
need to pass to the second part of this note, namely to test our 
definitions on some fractal sets.

In this paper we confine our attention to fractals in $\br$, namely 
to totally disconnected compact subsets of $\br$ with no isolated 
points. Even though fractals in $\br$ are not interesting from the 
point of view of fractal diffusions, they constitute a quite general 
class for our purposes, allowing very general situations from the 
point of view of Hausdorff-Besicovitch theory. Indeed this class, or 
better the class of their complements, the so called fractal strings, 
constituted the first playground  of the analysis of Lapidus and 
collaborators \cite{LavF}. 

The study of fractals in $\br^{n}$ from the point of view of 
noncommutative geometry will be the object of a forthcoming paper 
\cite{GuIs10}.

A simple spectral triple, which we call the ``lacunary'' spectral
triple, can be associated to fractals in $\br$, following the analysis
of Connes in \cite{Co}, IV.3.$\eps$.  We show here some measure
theoretic properties of such a triple.  In particular, for any fractal
in $\br$, its upper box dimension coincides with the Hausdorff
dimension $d$ of its lacunary spectral triple.  When the fractal is
$d$-Minkowski measurable, a result in \cite{LaPo} implies that the
singular traceability exponent is unique and equal to $D$, and the
corresponding functional is the Hausdorff functional.  By making use
of recent results of He and Lapidus \cite{HeLa}, we also prove that
for $h$-Minkowski measurable fractals the singular traceability
exponent is unique, the corresponding functional being not necessarily
Hausdorff.

However, the lacunary spectral triple does not reconstruct the
original metric on the fractal unless the fractal has zero Lebesgue
measure.  Following an idea of Connes \cite{Cpc}, we propose here a
new spectral triple for a wide class of fractals in $\br$, which we
call limit fractals, and can be seen as a subclass of the so called
random fractals \cite{random}.  Such a triple has the advantage of
retaining all the measure theoretic properties of the lacunary
spectral triple for limit fractals, and moreover reconstructs the
original metric without any further assumption.  This ``complete''
spectral triple can be described as a direct sum of the lacunary
spectral triple above, and a ``filled'' spectral triple, which we
proposed in \cite{GuIs8} as a spectral triple for limit fractals in
$\br^{n}$.  Indeed, all properties of limit fractals we prove here for
the complete spectral triple are valid for the direct summand
triples.

The simplest case of a limit fractal is a self-similar fractal.  We
show that for any self-similar fractal, the singular traceability
exponent is unique, and the associated Hausdorff-Besicovitch
functional is indeed the Hausdorff functional, and coincides (up to a
constant) with the Hausdorff measure.  The uniqueness result above in
the case of the lacunary triple is not implied by the analogous result
in the Minkowski measurable case, since not all self-similar fractals
are Minkowski measurable (cf.  \cite{Gatzuras}).  The fact that the
Hausdorff functional for the lacunary triple reconstructs the
Hausdorff measure has been shown in \cite{Co} IV.3.$\eps$ for
Cantor-like fractals, and is proved here in the general case (Theorem
\ref{Thm:self-similar}).

Then we show that, for limit fractals, the value of the singular trace
on the elements $f|D|^{-\a}$, $f$ being a continuous function, does
not depend on the generalized limit procedure (measurability in the
sense of Connes), namely the Hausdorff-Besicovitch functionals are
well defined.  We remark, in passing, that an analogous measurability
result has been recently proved by Kigami and Lapidus \cite{KiLa2}. 
They consider some class of self-similar fractals, for which a
Laplacian on the fractal can be constructed as a generator of a
Dirichlet form, and prove that the functional $f\to Tr_{\o}
(f\D^{-\a})$, where $\a>0$ is related to the self-similar dimension,
does not depend on $\o$.

Returning to this paper, we show that, in some cases (translation
fractals), the non commutative Hausdorff-Besicovitch functional
coincides with the restriction to the fractal of a
Hausdorff-Besicovitch measure on $\br$ (cf.  \cite{KaSa}).

For uniformly generated symmetric fractals, we are able to compute
explicitly $\subc$ and $\supc$, and this provides an evidence of a
classical interpretation for these numbers.  In fact in \cite{GuIs11}
we define upper and lower pointwise tangential dimensions for fractals
in $\br^{n}$.  These dimensions are defined by means of the box
dimensions of the balls of the tangent sets at a point, where a
tangent set of $F$ at $x$ is any limit, for $\l\to\infty$, of the
$\l$-dilations of $F$ around $x$, in a suitable topology.  We show in
\cite{GuIs11} that, for the uniformly generated symmetric fractals
considered here (cf.  Theorem \ref{Thm:delta*symm}), the upper, resp. 
lower, tangential dimension does not depend on $x$, and coincides with
$\supc$, resp.  $\subc$.  We notice here that while our motivation for
introducing the tangential dimensions was the attempt of finding a
classical counterpart of $\subc$ and $\supc$, the description given
above has been largely influenced by the notion of micro-set of
Furstenberg, as we heard it in his talk at Graz \cite{Fur}.

We conclude by mentioning that two of the results proved here have an
interest in the general study of singular traces.  The first is the
H\"older inequality, which we prove here for a general singular trace. 
In contrast with the Cauchy-Schwarz inequality, whose proof is purely
algebraic, H\"older inequality requires the characterization of
singular traces contained in \cite{GuIs1}.

The second is the complete description of the singular traceability
exponents, which is based, and indeed generalises, the
characterization of singular traceability given in \cite{AGPS}. 
Recently, we became aware of a related result contained in
\cite{DFWW}, where non-positive trace functionals on $B(H)$ are
studied.  Put together, these results suggest that the exponents of
singular traceability for positive singular traces should coincide
with those for non positive ones.

This paper is divided into four sections. The first two sections 
concern integration for spectral triples, the first containing the 
necessary technicalities on non increasing infinitesimal functions and 
the second the relevant results.

The last two sections concern fractals in $\br$, which are described 
in Section~3 from the classical point of view, while Section~4 
contains the results of our noncommutative analysis.  H\"older 
inequality for singular traces is proved in the Appendix.

The results of this work have been presented in several international 
conferences in the period 2000-2001. Some of them have been announced 
in the Proceedings of a Conference in Siena, July 2000, \cite{GuIs8}.

\section{Non-increasing infinitesimal functions.}\label{NonIncreasing}

 As is well known, we may associate, via non-increasing rearrangement, 
 a non-increasing infinitesimal function $\m$ defined on the interval 
 $[0,\infty)$ to any compact operator on a Hilbert space.

 In this section we treat some properties of the functions in this 
 class, which we will extensively use in the following sections in 
 order to get results concerning compact operators and singular 
 traces.

 Let $M$ be the class of non-increasing right-continuous infinitesimal
 functions $\m$ defined on the interval $[0,\infty)$, and $F$ be the
 class of non-decreasing right-continuous functions $f$ on $\br$,
 which are bounded from below and unbounded from above.  Clearly, the
 map
 \begin{equation}\label{f-mu}
	\mu(x)\mapsto f(t)=-\log\m(e^{t}),
 \end{equation}
 and $\m(0) := \lim_{x\to0^{+}} \m(x)$, gives a bijection between
 these two classes.

 Given $f\in F$, we consider the following asymptotic indices:
 \begin{align*}
	\subc(f)&=\left(\lim_{h\to\infty}\limsup_{t\to\infty}
	\frac{f(t+h)-f(t)}{h}\right)^{-1}\cr
	\supc(f)&=\left(\lim_{h\to\infty}\liminf_{t\to\infty}
	\frac{f(t+h)-f(t)}{h}\right)^{-1}\cr
	\subd(f)&=\left(\limsup_{t\to\infty}
	\frac{f(t)}{t}\right)^{-1}\cr
	\supd(f)&=\left(\liminf_{t\to\infty}
	\frac{f(t)}{t}\right)^{-1}	
 \end{align*}

 According to the previously mentioned correspondence between $M$ and 
 $F$, we shall write $\subd(\mu)$, $\supd(\mu)$, $\subc(\mu)$, 
 $\supc(\mu)$ as well.  Let us observe that these last two indices are 
 the Matuszewska indices of the function $\m$ (cf.  e.g. \cite{BGT}).  
 Some of the properties below may be known, but we prove them for the 
 sake of completeness.

 \subsection{Properties of the asymptotic indices}

 \begin{Prop}\label{constants}
	For any $f\in F$, the limits in the definitions of $\subc$, $\supc$ 
	exist, and the following relations hold:
	\begin{align*}
		\subc(f) &=\left(\inf_{h>0}\limsup_{t\to\infty} 
		\frac{f(t+h)-f(t)}{h}\right)^{-1} \leq\subd(f) \cr 
		&\leq\supd(f) \leq \left(\sup_{h>0}\liminf_{t\to\infty} 
		\frac{f(t+h)-f(t)}{h}\right)^{-1} = \supc(f).
	\end{align*}
 \end{Prop}

 \begin{proof}
	Let us set $\supg(h)=\limsup_{t\to\infty} \frac{f(t+h)-f(t)}{h}$.  
	Then we have
	\begin{align}\label{monotone}
	\supg(nh) 
	&= \limsup_{t\to\infty} \frac{1}{n} \sum_{k=1}^{n} 
	\frac{f(t+kh)-f(t+(k-1)h)}{h} \cr
	&\leq \frac{1}{n} \sum_{k=1}^{n} 
	\limsup_{t\to\infty} \frac{f(t+kh)-f(t+(k-1)h)}{h} = \supg(h).
	\end{align}
		
	Let us denote by $[s]_{r}$ the $r$-integer part of $s$, namely 
	$[s]_{r}=r[\frac{s}{r}]$.  Since $f$ is increasing, for any 
	infinite sequence $h_{n}$ we get
	$$
	\limsup_{n}\supg(h_{n})=\limsup_{n}\supg([h_{n}]_{r}),
	\quad
	\liminf_{n}\supg(h_{n})=\liminf_{n}\supg([h_{n}]_{r}).
	$$
	Now assume, ad absurdum, that $\limsup_{h\to\infty} \supg(h) > 
	\liminf_{h\to\infty} \supg(h)$.  Then we can find two sequences 
	$h_{i}$ and $h'_{i}$ such that $\supg(h_{i}) \to 
	\liminf_{h\to\infty}\supg(h)$, $\supg(h'_{i}) \to 
	\limsup_{h\to\infty}\supg(h)$, and $\supg(h_{i}) < \supg(h'_{j})$ 
	for any $i,j$. But
	$$
	\lim_{j}\supg(h'_{j}) = \lim_{j}\supg([h'_{j}]_{h_{i}}) = 
	\lim_{j}\supg(h_{i}\left[\frac{h'_{j}}{h_{i}}\right]) \leq 
	\supg(h_{i}),
	$$
	which is absurd.  The proof for $\subg$ is analogous.  \\
	Now we prove the relations.  By equation (\ref{monotone}) and the 
	existence of the limit of $\supg(h)$ when $h\to\infty$, we have 
	that $\lim_{h\to\infty}\supg(h)=\inf_{h\geq 0}\supg(h)$, namely 
	the first equation.  \\
	Since $f$ is increasing, we	get 
	$$
	\limsup_{t\to\infty}\frac{f(t)}{t} = 
	\limsup_{t\to\infty}\frac{f([t]_{h})}{[t]_{h}},
	$$
	therefore, setting $a_{k}=\frac{f(kh)-f((k-1)h)}{h}$, we have, for 
	any $h>0$,
	\begin{align*}
		\limsup_{t\to\infty}\frac{f(t)}{t} 
		& =\limsup_{n}\frac{f(nh)}{nh} = \limsup_{n} \left( 
		\frac{1}{n} \sum_{k=1}^{n} a_{k} + \frac{f(0)}{nh} \right) \cr 
		&\leq \limsup_{k}a_{k}\leq\limsup_{t\to\infty} 
		\frac{f(t+h)-f(t)}{h},
	\end{align*}
	or, equivalently, 
	$$
	\limsup_{t\to\infty}\frac{f(t)}{t}\leq \inf_{h\geq 0}\supg(h),
	$$
	which is the first inequality of the statement.
	\\
	The last equality, resp. inequality of the statement are proved analogously.
 \end{proof}

 \begin{Lemma}\label{Lemma:doublelim}
    \begin{align}
		\subc(f)^{-1}&=\limsup_{t,h\to\infty} \frac{f(t+h)-f(t)}{h},\\
		\supc(f)^{-1}&=\liminf_{t,h\to\infty} \frac{f(t+h)-f(t)}{h}.
    \end{align}
 \end{Lemma}
 \begin{proof}
     Setting $\f(t,h) := \frac{f(t+h)-f(t)}{h}$, we have to show
     $\limsup_{t,h\to\infty} \f(t,h) = \lim_{h\to\infty}
     \limsup_{t\to\infty} \f(t,h)$.  Assume
     $\lim_{h\to\infty}\limsup_{t\to\infty}\f(t,h)=L\in\br$.  Let
     $\eps>0$, then there is $h_{\eps}>0$ such that, for any
     $h>h_{\eps}$, $\limsup_{t\to\infty}\f(t,h) > L-\eps/2$, hence,
     for any $t_{0}>0$ there is $t=t(h,t_{0})$, such that
     $\f(t,h)>L-\eps$.  Hence, for any $h_{0}>0$, $t_{0}>0$ there
     exist $h>h_{0}$, $t>t_{0}$ such that $\f(t,h)>L-\eps$, namely
     $\limsup_{t,h\to\infty}\f(t,h)\geq L$.  Conversely, assume
     $\limsup_{t,h\to\infty}\f(t,h)=L'\in\br$, and choose $t_{n}$,
     $h_{n}$ such that $\lim_{n\to\infty}\f(t_{n},h_{n})=L'$.  For any
     $r>0$, let us denote by $\{s\}_{r}:= r
     \left\{\frac{s}{r}\right\}$, where $\{s\}$ is the least integer
     no less than $s$.  Then, for any $h>0$, with $p$ denoting
     $\frac{\{h_{n}\}_{h}}{h}$, we have
     \begin{align*}
	\f(t_{n},h_{n}) & \leq \frac{\{h_{n}\}_{h}}{h_{n}}
	\f(t_{n},\{h_{n}\}_{h}) = \frac{\{h_{n}\}_{h}}{h_{n}}\frac{1}{p}
	\sum_{j=0}^{p-1}\f(t_{n}+jh,h)\\
	&\leq \frac{\{h_{n}\}_{h}}{h_{n}}\max_{j=0\dots p-1}
	\f(t_{n}+jh,h) \leq \frac{\{h_{n}\}_{h}}{h_{n}}\sup_{t\geq
	t_{n}}\f(t,h).
    \end{align*}
    Hence, for $n\to\infty$, we get
    $L'\leq\limsup_{t\to\infty}\f(t,h)$, which implies the equality. 
    The other cases are treated analogously.
 \end{proof}

 Now we introduce the notion of eccentricity for a function $\m\in M$.  
 It is motivated by the fact that a positive compact operator (cf.  
 next section) is singularly traceable if and only if its eigenvalue 
 function is eccentric.

 \begin{Dfn}\label{Dfn:ecc}
	 Given a function $\m\in M$, we define its integral function $S$ 
	 as
	 \begin{equation*}
		 S(x) = 
		 \begin{cases}
			 \sinc(x):=\int_0^x \m(y)dy&\m\not\in L^1[0,\infty)\cr
			 \sdec(x):=\int_x^\infty\m(y)dy&\m\in L^1[0,\infty).
		 \end{cases}
	 \end{equation*}
	 A function $\m\in M$ is eccentric if $1$ is a limit point, when 
	 $x\to\infty$, of the function
	 $$
	 \frac{S(\l x)}{S(x)},
	 $$
	 for some $\l>1$.  Note that if it is true for one $\l$, it is 
	 indeed true for any $\l>1$, cf.  \cite{GuIs1}.
 \end{Dfn}

 It was proved in \cite{GuIs5} that $\supd(\m)=1$ is a sufficient 
 condition for $\m$ to be eccentric.

 \begin{Thm}\label{Thm:dimension} 
	 Let $\m$ be a function in $M$.  Then the following (possibly 
	 infinite) quantities coincide with $\supd(\m)$:
	 \begin{align*}
		 d_1 & = \sup \{\a>0: \limsup_{x\to\infty} 
		 \frac{\int_{0}^x\m(t)^\a dt}{\log x} = \infty \} \\
		 d_2 &=\inf\{\a>0:\lim_{x\to\infty}\frac{\int_{0}^x\m(t)^\a 
		 dt}{\log x}=0\} \\
		 d_3 &=\inf\{\a>0:\int_{0}^x\m(t)^\a dt<\infty\} \\
		 d_4 &=\left(\liminf_{t\to\infty} 
		 \frac{\log1/\m(t)}{\log t}\right)^{-1}
	 \end{align*}
 \end{Thm}

 \begin{proof} 
	 Set
	 \begin{align*}
		 \O_1 &=\{\a>0:\limsup_{x\to\infty}\frac{\int_{0}^x\m(t)^\a 
		 dt}{\log x}=\infty\} \\
		 \O_2 &=\{\a>0:\lim_{x\to\infty}\frac{\int_{0}^x\m(t)^\a 
		 dt}{\log x}=0\} \\
		 \O_3 &=\{\a>0:\int_{0}^x\m(t)^\a dt<\infty\}.
	 \end{align*}
	 \\
	 $d_1\leq d_2$.  If $\a\in\O_1$, then $(0,\a]\supseteq\O_1$, if 
	 $\b\in\O_2$, then $[\b,\infty)\supseteq\O_2$, and 
	 $\O_1\cap\O_2=\emptyset$, hence $\O_1$ and $\O_2$ are separated 
	 classes: $\O_1\leq\O_2$.  \\
	 $d_2\leq d_3$. If $\a\in\O_3$ then $\a\in \O_2$, namely $\O_3\subseteq\O_2$.
	 \\
	 $d_3\leq d_4$.  Let $a(t)=\frac{\log1/\m(t)}{\log(t)}$, namely 
	 $\m(t)=t^{-a(t)}$ and $\liminf_{t\to\infty} a(t)=1/d_4$.  If 
	 $\a>d_4$, then $\liminf_{t\to\infty}\a a(t)=\a/d_4>1$, hence 
	 there exists $\b>1$ such that $\a a(t)\geq\b$ for $t$ 
	 sufficiently large.  Therefore
	 $$
	 \int_{0}^\infty\m(t)^\a dt=\int_{0}^\infty t^{-\a a(t)}dt\leq const+
	 \int_{0}^\infty t^{-\b}dt<\infty,
	 $$
	 which implies $\a\in\O_3$, namely $(d_4,\infty)\subset\O_3$.
	 \\
	 $d_4\leq d_1$.  We may assume $d_4>0$, namely $1/d_4<\infty$.  
	 Now let $\ell\in[0,\infty)$ be a limit point, for $t\to\infty$, of the 
	 function $\frac {\log1/\m(t)} {\log(t)}$, namely $\ell_k:=\frac 
	 {\log1/\m(t_k)} {\log {t_k}} \to \ell$, for a suitable sequence 
	 $t_k\to\infty$.  We have $\m(t_k) = t_k^{-\ell_k}$.  Let now 
	 $\a<1/\ell$ (possibly $1/\ell=\infty$), namely $\a \ell_k\to\a \ell<1$, and 
	 choose $\eps>0$ such that $\a \ell_k\leq 1-\eps$ eventually.  Then
	 $$
	 \int_{0}^{t_k}\m(t)^\a dt \geq t_k\m(t_k)^\a = t_k\cdot t_k^{-\a 
	 \ell_k} \geq t_k^\eps.
	 $$
	 Therefore 
	 $$
	 \frac{\int_{1}^{t_k}\m(t)^\a}{\log t_k} \geq \frac{t_k^\eps}{\log 
	 t_k} \to\infty,
	 $$
	 which means  that $\a\in \O_1$, i.e. $(0,1/\ell)\subseteq \O_1$ and 
	 $d_1\geq 1/\ell$.
	 \\
	 The equality $d_4=\supd$ follows immediately from the definitions.
 \end{proof}

 We want to show that $\m^{\a}$ is eccentric {\it iff} 
 $\a\in[\subc,\supc]\cap(0,\infty)$, thus giving a meaning to the 
 quantities $\subc,\supc$.
 
 \begin{Thm}\label{Thm:subcsupc} 
	 Let $\m$ be in $M$, $\g$ be a positive number.  If $\m^{\g}$ is 
	 eccentric, then $\g\in[\subc,\supc]$.
 \end{Thm}

 \begin{proof}
	 Let $\a<\subc(\m)$.  The first equality in Proposition 
	 \ref{constants} says that $\a<\subc(\m)$ if and only if 
	 $\a<\sup_{h>0} \left( \limsup_{t\to\infty} \frac{f(t+h)-f(t)}{h} 
	 \right)^{-1}$, namely if and only if there exists $h>0$ such that 
	 $\a < \left( \limsup_{t\to\infty} \frac{f(t+h)-f(t)}{h} 
	 \right)^{-1}$, or, equivalently, using the function $\m$ 
	 associated with $f$, if there exists $\l>1$ for which
	 $$
	 \l\liminf_{t\to\infty}\frac{\m(\l t)^{\a}}{\m(t)^{\a}}>1.
	 $$
	 Now observe that, by the inequalities in Proposition 
	 \ref{constants}, $\a<\supd$, hence, by Theorem 
	 \ref{Thm:dimension}, $\m^{\a}\not\in L^{1}$.  Therefore
	 $$
	 \frac{S_{\m^{\a}}(\l x)}{S_{\m^{\a}}(x)} 
	 =\frac{\l\int_{0}^{x}\m(\l t)^{\a} dt}{\int_{0}^{x}\m(t)^{\a} dt} 
	 =\l\frac{\int_{0}^{x}\left(\frac{\m(\l t)}{\m(t)}\right)^{\a} 
	 \m(t)^{\a} dt}{\int_{0}^{x}\m(t)^{\a} dt},
	 $$
	 hence
	 $$
	 \liminf_{x\to\infty}\frac{S_{\m^{\a}}(\l x)}{S_{\m^{\a}}(x)} \geq 
	 \l\liminf_{x\to\infty}\left(\frac{\m(\l x)}{\m(x)}\right)^{\a}>1,
	 $$
	 which implies that $\m^{\a}$ is not eccentric.
	 \\
	 The proof for $\a>\supc(\m)$ is analogous.
 \end{proof}

 To prove the converse direction, we need some preliminary results.

 \begin{Prop} 
	 1 is a limit point of $\frac{S(x)}{S(2x)}$ if and only if 
	 $\liminf\frac{x\m(x)}{S(x)}=0$.
 \end{Prop}

 \begin{proof}
	 Assume first $\m$ is not summable, i.e. $S(x)=\sinc(x)$.  Then the 
	 thesis follows by the following inequalities:
	 \begin{align*}
		 \frac{\sinc(2x)}{\sinc(x)}-1&\leq\frac{x\m(x)}{\sinc(x)}\\ 
		 \frac{2x\m(2x)}{\sinc(2x)} &\leq2\frac{\sinc(2x)-\sinc(x)}{\sinc(2x)} 
		 =2\left(1-\frac{\sinc(x)}{\sinc(2x)}\right).
	 \end{align*}
	 When $\m$ is summable, i.e. $S(x)=\sdec(x)$, we have, analogously,
	 \begin{align*}
		 1-\frac{\sdec(2x)}{\sdec(x)}&\leq\frac{x\m(x)}{\sdec(x)}\\ 
		 \frac{2x\m(2x)}{\sdec(2x)} &\leq 
		 2\left(\frac{\sdec(x)}{\sdec(2x)}-1\right)
	 \end{align*}
	 and the thesis follows.
 \end{proof}

 \begin{Prop}\label{Prop:lessthansubc} 
	 If $\inf_{t}\frac{t\m(t)}{\sinc(t)}= k>0$ for any $t>0$, then 
	 $\subc>1$.
 \end{Prop}

 \begin{proof}
    Since $\m$ is the derivative of $\sinc$, the hypothesis means that
    $$
    \frac{d}{dt} \log\sinc(t)\geq\frac{k}{t}, \forall t.
    $$
    Integrating on the interval $[x,\l x]$ one gets
    $$
    \sinc(\l x)\geq \l^{k}\sinc(x).
    $$
    Since $\frac{x\m(x)}{\sinc(x)}\leq 1$, one obtains
    \begin{equation*}
	\frac{\m(x)}{\m(\l x)}\leq\frac{\frac{\sinc(x)}{x}}{\frac{k\sinc(\l 
	x)}{\l x}} \leq \frac{\l^{1-k}}{k}.
    \end{equation*}
    As a consequence,
    \begin{align*}
	\subc^{-1}&=\lim_{\l\to\infty}\frac{1}{\log 
	\l}\limsup_{x\to\infty}\log\frac{\m(x)}{\m(\l x)}\cr 
	&\leq\lim_{\l\to\infty}\frac{-\log 
	k+(1-k)\log\l}{\log\l}=1-k<1
    \end{align*}
 \end{proof}

 \begin{Prop}\label{Prop:morethansupc} 
	 If $\inf_{t}\frac{t\m(t)}{\sdec(t)}= k>0$ for any $t>0$, then 
	 $\supc < 1$.
 \end{Prop}

 \begin{proof}
    Since $-\m$ is the derivative of $\sdec$, we may prove, in analogy 
    with the previous Proposition, that
    $$
    \sdec(\l x)\leq \l^{-k}\sdec(x).
    $$
    Since $\sdec(t)\geq\int_{t}^{2t}\m(s)ds\geq t\m(2t)$, which implies
    $\frac{x\m(2x)}{\sdec(x)}\leq 1$, one obtains
    \begin{equation*}
		\frac{\m(2\l x)}{\m(x)}\leq\frac{\frac{\sdec(\l x)}{\l 
		x}}{\frac{k\sdec(x)}{x}} \leq \frac{\l^{-1-k}}{k},
    \end{equation*}
    namely
    \begin{equation*}
		\frac{\m(x)}{\m(\l x)} \geq k\left(\frac{\l}{2}\right)^{1+k},
    \end{equation*}
	As a consequence,
	\begin{align*}
		\supc^{-1}&=\lim_{\l\to\infty}\frac{1}{\log 
		\l}\liminf_{x\to\infty}\log\frac{\m(x)}{\m(\l x)}\cr 
		&\geq\lim_{\l\to\infty}\frac{\log 
		k+(1+k)\log2+(1+k)\log\l}{\log\l}=1+k>1,
	\end{align*}
 \end{proof}

 \begin{Thm}\label{Thm:eccint} 
	 The set of the eccentricity exponents of $\m$ is the interval 
	 whose endpoints are $\subc(\m)$, $\supc(\m)$ and is relatively 
	 closed in $(0,\infty)$.
 \end{Thm}

 \begin{proof} 
	Assume $\m^{\g}$ is not eccentric.  \\
    If $\m^{\g}\not\in L^{1}$, then 
    $\inf_{t}\frac{t\m(t)^{\g}}{\int_{0}^{t}\m(s)^{\g}ds}>0$ for any 
    $t>0$.  Therefore, by Proposition \ref{Prop:lessthansubc}, 
    $1<\subc(\m^{\g})=\frac{1}{\g}\subc(\m)$, i.e. $\g<\subc(\m)$.  \\
    If $\m^{\g}\in L^{1}$, then 
    $\inf_{t}\frac{t\m(t)^{\g}}{\int_{t}^{\infty}\m(s)^{\g}ds}>0$ for 
    any $t>0$.  Therefore, by Proposition \ref{Prop:morethansupc}, 
    $1>\supc(\m^{\g})=\frac{1}{\g}\supc(\m)$, i.e. $\g>\supc(\m)$.  \\
    The converse implication is contained in Theorem \ref{Thm:subcsupc}.
 \end{proof}

 \subsection{Direct sums of infinitesimal functions.}

 Let $f$, $g$ be two real valued measurable functions defined on the
 measure spaces $A$ and $B$ respectively.  We say (cf.  \cite{BS})
 that $g$ is a rearrangement of $f$ if there is a measure preserving
 bijection $\f$ from the support of $f$ to the support of $g$ and
 $f=g\circ\f$.  The nonincreasing rearrangement $f^{*}$ of $f$ is
 defined as the unique non-increasing, right continuous rearrangement
 of $f$ on $[0,\infty)$ with the Lebesgue measure.

 It is known that $f^{*}$ can be defined as $f^{*}(t):=\inf\{s\geq0:
 \l_f(s)\leq t\}$, $t\geq0$, where $\l_f(t)$ is the measure of
 $\{x:f(x)>t\}$.

 Consider now the following binary operation on $M$: let $\a$, $\b\in
 M$, and set $\a\oplus\b$ to be the nonincreasing rearrangement of
 $\tilde\a+\tilde\b$, where $\tilde\a$ and $\tilde\b$ have disjoint
 supports and $\a$, resp.  $\b$ is the nonincreasing rearrangement of
 $\tilde\a$, resp.  $\tilde\b$.  This operation is well defined,
 namely does not depend on the rearrangements $\tilde\a$ and
 $\tilde\b$.  Indeed $\l_{\tilde\a+\tilde\b} = \l_{\tilde\a} +
 \l_{\tilde\b} = \l_{\a}+\l_{\b}$.

 The need for this operation relies on Proposition \ref{mu-oplus}
 below.

 \begin{Prop}\label{sumeqsup}
	Let $\a$, $\b$ be elements of $M$. Then
	\begin{align}
		\subc(\a\oplus\b)&=\subc(\a\vee\b)\label{a1}\\
		\supc(\a\oplus\b)&=\supc(\a\vee\b)\label{a2}\\
		\subd(\a\oplus\b)&=\subd(\a\vee\b)\label{a3}\\
		\supd(\a\oplus\b)&=\supd(\a\vee\b)\label{a4}.
	\end{align}
 \end{Prop}

 \begin{proof}
	We have $\a=(\tilde\a)^{*}\leq(\tilde\a+\tilde\b)^{*}=\a\oplus\b$ and 
	analogously for $\b$, therefore we get 
	\begin{equation}\label{a5}
	(\a\vee\b)(x)\leq(\a\oplus\b)(x).
	\end{equation}
	Moreover, since $(a+b)^{*}(s+t)\leq a^{*}(s)+b^{*}(t)$ (cf.  
	\cite{BS}), we have
	\begin{equation}\label{a6}
	(\a\oplus\b)(2x)\leq 2(\a\vee\b)(x),
	\end{equation}
	from which equalities (\ref{a1}) to (\ref{a4}) follow immediately.
 \end{proof}

 \begin{Prop}\label{estimates}
	Let $\a$, $\b$ be elements of $M$. Then
	\begin{align}
		\subc(\a\oplus\b)&\geq\subc(\a)\wedge\subc(\b)\label{a7}\\
		\supc(\a\oplus\b)&\leq\supc(\a)\vee\supc(\b)\label{a8}\\
		\subd(\a\oplus\b)&\geq\subd(\a)\vee\subd(\b)\label{a9}\\
		\supd(\a\oplus\b)&=   \supd(\a)\vee\supd(\b)\label{a10}.
	\end{align}
 \end{Prop}

 \begin{proof}	
	We have
	\begin{align*}
	\supd(\a\oplus\b)^{-1}
	&=\supd(\a\vee\b)^{-1}\cr
	&=\liminf_{x\to\infty}\frac{-\log(\a(x)\vee\b(x))}{\log x}\cr
	&=\liminf_{x\to\infty}\left(\frac{-\log\a(x)}{\log x}
	\wedge\frac{-\log\b(x)}{\log x}\right)\cr
	&=\supd(\a)^{-1}\wedge\supd(\b)^{-1},
	\end{align*}
	which shows equation (\ref{a10}). Inequality (\ref{a9}) is proved 
	analogously.
	\\
	Using (\ref{a1}), we get
	\begin{multline}\label{eq:oplus}
		\subc(\a\oplus\b)^{-1} 
		=\lim_{\l\to\infty}\frac{1}{\log \l}\limsup_{x\to\infty}
		\log\frac{(\a\vee\b)(x)}{(\a\vee\b)(\l x)}\cr 
		=\lim_{\l\to\infty}\limsup_{x\to\infty}
		\frac{(\log\a(x)\vee\log\b(x))+(-\log\a(\l x)\wedge-\log\b(\l x))}
		{\log \l},
	\end{multline}
	from which, observing that $(a\vee b)+(c\wedge d)\leq (a+c)\vee 
	(b+d)$, we get (\ref{a7}).  Making use of 
	$(a+c)\wedge(b+d)\leq(a\vee b)+(c\wedge d)$, inequality (\ref{a8}) 
	is proved.
 \end{proof}	
	
 \begin{Thm}\label{c(oplus)}
	The interval $[\subc(\a\oplus\b),\supc(\a\oplus\b)]$ is contained 
	in the interval 
	$[\subc(\a)\wedge\subc(\b),\supc(\a)\vee\supc(\b)]$.  This 
	estimate is optimal, namely the equality may happen, but the 
	interval may shrink to a point $(=\supd(\a\oplus\b))$ in some cases.
 \end{Thm}

 \begin{proof} 
	The first statement immediately follows from Proposition 
	\ref{estimates}.  The estimate is optimal e.g. in the case of 
	Proposition \ref{dim.tang.eq}. The interval shrinks to a point in the 
	following example.
	\\
	Making use of the identification (\ref{f-mu}), the two functions may be 
	equivalently described by two non-decreasing functions $f$, $g$. 
	By equations (\ref{a1}), (\ref{a2}), the direct sum now corresponds to 
	$f\wedge g$.  Let us choose $f$ and $g$ as follows: choose a 
	sequence of intervals $I_{n}=(a_{n},a_{n+1}]$, $n\in\bn$, with 
	increasing length, set $f(t)=t$ for $t\in I_{n}$, $n$ even, and 
	$f(t)=a_{n+1}$ for $t\in I_{n}$, $n$ odd.  It is easy to check that 
	$\subc(f)=0$ and $\supc(f)=+\infty$.  Conversely, choose $g(t)=t$ for 
	$t\in I_{n}$, $n$ odd, and $g(t)=a_{n+1}$ for $t\in I_{n}$, $n$ even.  
	Again, $\subc(g)=0$ and $\supc(g)=+\infty$.  Moreover, $(f\wedge g)(t)=t$, 
	therefore the corresponding interval shrinks to the point $\{1\}$.
 \end{proof}

 For the application to fractals in Section 4, we need a refinement of 
 the previous Theorem.
 
 \begin{Prop}\label{dim.tang.eq}
     Let $\a$, $\b$ be elements of $M$, and assume there are $A,B\geq
     1$ such that $\frac{1}{A}\a(Bx) \leq \b(x) \leq A\a(\frac{x}{B})$, 
     for all $x$ large enough. Then 
     	\begin{align}
		\subc(\a\oplus\b)&=\subc(\a) =\subc(\b)\\
		\supc(\a\oplus\b)&=\supc(\a) = \supc(\b).
	\end{align}
 \end{Prop}
 \begin{proof}
     From the hypotheses we get, for all $x$ large enough,
     $$
	 -\log A + \log \a(Bx) \leq \log\a(x)\vee\log\b(x) \leq \log A + \log 
	 \a(\frac{x}{B}).
     $$
     Therefore, from equation (\ref{eq:oplus}), we get
     \begin{align*}
	\subc(\a\oplus\b)^{-1}
	&=\lim_{\l\to\infty}\limsup_{x\to\infty}
	\frac{(\log\a(x)\vee\log\b(x))+(-\log\a(\l x)\wedge-\log\b(\l
	x))} {\log \l} \\
	& \leq \lim_{\l\to\infty}\limsup_{x\to\infty} \frac{2\log A +
	\log \a(\frac{x}{B}) - \log\a(\l Bx)} {\log \l} = \subc(\a)^{-1},
    \end{align*}
    and in like manner we obtain the reversed inequality. All the other 
    statements are proved in the same way.
 \end{proof}

 \section{Integration for spectral triples}\label{NonCommInt}

 In this section we apply the results of the previous section in order
 to study traceability properties of compact operators and then to
 interpret them in the framework of Alain Connes' Noncommutative
 Geometry.

 \subsection{Singular traceability}

 The theory of singular traces on $\cb(\ch)$, namely positive trace 
 functionals vanishing on the finite rank projections, was developed 
 by Dixmier \cite{Dixmier}, who first showed their existence, and then 
 in \cite{Varga}, \cite{AGPS}.  For the theory of non-positive traces 
 see \cite{DFWW}.  For generalizations to von Neumann algebras and 
 C$^{*}$-algebras see \cite{GuIs1,GuIs4,DPSS,BeFa,CPS}.

 Any tracial weight is finite on an ideal contained in $\ck(\ch)$ and 
 may be decomposed as a sum of a singular trace and a multiple of the 
 normal trace.  Therefore the study of (non-normal) traces on 
 $\cb(\ch)$ is the same as the study of singular traces.

 Moreover, because of singularity and unitary invariance, a 
 singular trace depends only on the eigenvalue asymptotics, 
 namely, if $a$ and $b$ are positive compact operators on $\ch$ and 
 $\m_n(a)=\m_n(b)+o(\m_n(b))$, $\m_n$ denoting the $n$-th eigenvalue, 
 then $\t_{\o}(a)=\t_{\o}(b)$ for any singular trace $\t_{\o}$.

 The main problem about singular traces is therefore to detect which 
 asymptotics may be ``summed'' by a suitable singular trace, that is 
 to say, which operators are singularly traceable.

 In order to state the most general result in this respect we need 
 some notation.

 Let $a$ be a compact operator.  Then we denote by 
 $\{\m_n(a)\}_{n=0}^{\infty}$ the sequence of the eigenvalues of 
 $|a|$, arranged in non-increasing order and counted with 
 multiplicity, and by $\m_{a}$ the corresponding eigenvalue function, 
 which is equal to $\m_{k}(a)$ on the interval $[k,k+1)$ for any $k$.  
 We denote the corresponding integral function $S_{\m_{a}}$, defined 
 in the previous section, simply by $S_{a}$.
 
 A compact operator is called {\it singularly traceable} if there 
 exists a singular trace which is finite non-zero on $|a|$.  We 
 observe that the domain of such singular trace should necessarily 
 contain the ideal $\cai(a)$ generated by $a$.  Then the following 
 theorem holds.

 \begin{Thm}\label{eccsingtrac}\cite{AGPS}
     A positive compact operator $a$ is singularly traceable $iff$
     $\m_{a}$ is eccentric (cf.  Definition \ref{Dfn:ecc}).  In this
     case there exists a sequence $x_k\to\infty$ such that, for any
     generalised limit $\Lim_{\omega}$ on $\ell^\infty$, the positive
     functional
     $$
     \tau_\omega (b) = 
     \begin{cases} 
	 \Lim_{\omega}\left(\left\{\frac{S_{b}(x_k)}{S_{a}(x_k)}\right\}
	 \right) 
	 &\quad b \in \cai(a)_+ \\
	 +\infty&\quad b \not\in \cai(a),\ b>0,
     \end{cases}
     $$
     is a singular trace whose domain is the ideal $\cai(a)$ generated
     by $a$.
 \end{Thm}

 The best known eigenvalue asymptotics giving rise to a singular trace 
 is $\m_{n}\sim\frac1n$, which implies $S(x)\sim\log x$.  The corresponding 
 logarithmic singular trace is generally called Dixmier trace.

 \begin{Dfn} 
	 If $a\in\cpt$ we define
	 $\subc(a)=\subc(\m_{a})$, $\supc(a)=\supc(\m_{a})$, 
	 $\subd(a)=\subd(\m_{a})$, $\supd(a)=\supd(\m_{a})$.
	 We say that $\a>0$ is an exponent of singular traceability for $a$ 
	 if $|a|^{\a}$ is singularly traceable.
 \end{Dfn}
 
 \begin{Thm}\label{Thm:SingTracExp} 
     Let $a$ be a compact operator.  Then, the set of singular
     traceability exponents is the closed interval in $(0,\infty)$
     whose endpoints are $\subc(a)$ and $\supc(a)$.  In particular, if
     $\supd(a)$ is finite nonzero, it is an exponent of singular
     traceability.
 \end{Thm}

 \begin{proof}
	The statement follows by Theorems \ref{Thm:eccint}, \ref{eccsingtrac}.
 \end{proof}

 Note that the interval of singular traceability may be $(0,\infty)$,
 as shown in \cite{GuIs8}.
 
 In \cite{GuIs12} the previous Theorem has been generalised to any 
 semifinite factor, and some questions concerning the 
 domain of a singular trace have been considered.

 \subsection{Singular traces and spectral triples}

 In this section we shall discuss some notions of dimension in 
 noncommutative geometry in the spirit of Hausdorff-Besicovitch 
 theory.

 As is known, the measure for a noncommutative manifold is defined 
 via a singular trace applied to a suitable power of some geometric 
 operator (e.g. the Dirac operator of the spectral triple of Alain 
 Connes).  Connes showed that such procedure recovers the usual volume 
 in the case of compact Riemannian manifolds, and more generally the 
 Hausdorff measure in some interesting examples \cite{Co}, Section 
 IV.3.

 Let us recall that $(\ca,\ch,D)$ is called a {\it spectral triple} when 
 $\ca$ is an algebra acting on the Hilbert space $\ch$, $D$ is a self 
 adjoint operator on the same Hilbert space such that $[D,a]$ is 
 bounded for any $a\in\ca$, and $D$ has compact resolvent.  In the 
 following we shall assume that $0$ is not an eigenvalue of $D$, the 
 general case being recovered by replacing $D$ with 
 $D|_{\ker(D)^\perp}$.  Such a triple is called $d^+$-summable, $d\in 
 (0,\infty)$, when $|D|^{-d}$ belongs to the Macaev ideal 
 $\cl^{1,\infty}=\{a:\frac{S_{a}^{\uparrow}(t)}{\log t}<\infty\}$.  
 
 The noncommutative version of the integral on functions is given by 
 the formula $\Tr_\omega(a|D|^{-d})$, where $\Tr_\omega$ is the 
 Dixmier trace, i.e. a singular trace summing logarithmic divergences.  
 By the arguments below, such integral can be non-trivial only if $d$ 
 is the Hausdorff dimension of the spectral triple, but even this 
 choice does not guarantee non-triviality.  However, if $d$ is finite 
 non-zero, we may always find a singular trace giving rise to a 
 non-trivial integral.

 \begin{Thm} \label{Thm:trace}
	 Let $(\ca,\ch,D)$ be a spectral triple.  If $s$ is an exponent of 
	 singular traceability for $|D|^{-1}$, namely there is a singular 
	 trace $\t_{\o}$ which is non-trivial on the ideal generated by 
	 $|D|^{-s}$, then the functional $a\mapsto\t_{\o}(a |D|^{-s})$ is a 
	 trace state (Hausdorff-Besicovitch functional) on the algebra 
	 $\ca$. 
 \end{Thm}
 
 \begin{proof}
 	It is the same as the proof of Theorem 1.3 in \cite{CiGS1}, by 
 	making use of the H\"older inequality for singular traces proved in 
 	the Appendix.
 \end{proof}
 
 \begin{rem}
	When $(\ca,\ch,D)$ is associated to an $n$-dimensional compact 
	manifold $M$, or to the fractal sets considered in \cite{Co}, the 
	singular trace is the Dixmier trace, and the associated functional 
	corresponds to the Hausdorff measure.  This fact, together with the 
	previous theorem, motivates the following definition.
 \end{rem}
 
 \begin{Dfn}\label{Dfn:dimensions} 
	   Let $(\ca,\ch,D)$ be a spectral triple, $\Tr_{\o}$ the Dixmier 
	   trace.  \itm{i} We call $\a$-dimensional Hausdorff functional 
	   the map $a\mapsto Tr_\omega(a |D|^{-\a})$; \itm{ii} we call 
	   (Hausdorff) dimension of the spectral triple the number
	   $$
	   d(\ca,\ch,D) = \inf \{ d>0: |D|^{-d} \in \cl^{1,\infty}_{0} \} = 
	   \sup \{ d>0: |D|^{-d} \not\in \cl^{1,\infty} \},
	   $$
	   where $\cl^{1,\infty}_{0} = \{a : 
	   \frac{S_{a}^{\uparrow}(t)}{\log t} \to 0 \}$.  \itm{iii} we 
	   call minimal, resp.  maximal dimension of the spectral triple 
	   the quantity $\subc(|D|^{-1})$, resp.  $\supc(|D|^{-1})$.  
	   \itm{iv} For any $s$ between the minimal and the maximal 
	   dimension, we call the corresponding trace state on the algebra 
	   $\ca$ a Hausdorff-Besicovitch functional on $(\ca,\ch,D)$.
 \end{Dfn}
  
 \begin{Thm}\label{unisingtrac}
     \itm{i} $d(\ca,\ch,D)=\supd(|D|^{-1})$. 
     \itm{ii} $d := d(\ca,\ch,D)$ is the unique exponent, if any, such
     that $\ch_{d}$ is non-trivial.  
     \itm{iii} If $d\in (0,\infty)$, it is an exponent of singular
     traceability.
 \end{Thm}
 
 \begin{proof}
     $(i)$ The equality directly follows from Theorem 
     \ref{Thm:dimension}. \\
     $(ii)$ It follows easily from the definition. \\
     $(iii)$ It is a direct consequence of $(i)$
     and of Theorem \ref{Thm:eccint}.
 \end{proof}
 
 Let us observe that the $\a$-dimensional Hausdorff functional depends 
 on the generalized limit procedure $\omega$, however its value is 
 uniquely determined on the operators $a\in\ca$ such that $a|D|^{-d}$ 
 is measurable in the sense of Connes \cite{Co}. By an abuse of 
 language we call measurable such operators.  
 
 As in the commutative case, the dimension is the supremum of the 
 $\a$'s such that the $\a$-dimensional Hausdorff measure is everywhere 
 infinite and the infimum of the $\a$'s such that the $\a$-dimensional 
 Hausdorff measure is identically zero.  Concerning the non-triviality 
 of the $d$-dimensional Hausdorff functional, we have the same 
 situation as in the classical case.  Indeed, according to the 
 previous result, a non-trivial Hausdorff functional is unique (on 
 measurable operators) but does not necessarily exist.  In fact, if 
 the eigenvalue asymptotics of $D$ is e.g. $n\log n$, the Hausdorff 
 dimension is one, but the 1-dimensional Hausdorff measure gives the 
 null functional.
 
 However, if we consider all singular traces, not only the logarithmic 
 ones, and the corresponding trace functionals on $\ca$, as we said, 
 there exists a non trivial trace functional associated with 
 $d(\ca,\ch,D)\in(0,\infty)$, but $d(\ca,\ch,D)$ is not characterized 
 by this property.  In fact this is true if and only if the minimal 
 and the maximal dimension coincide.  A sufficient condition is the 
 following.
 
 \begin{Prop}\label{Prop:unique} 
	 Let $(\ca,\ch,D)$ be a spectral triple with finite non-zero 
	 dimension $d$.  If there exists $\lim \frac{\m_n(D^{-1})} 
	 {\m_{2n}(D^{-1})} \in(1,\infty)$, $d$ is the unique exponent of 
	 singular traceability of $D^{-1}$.
 \end{Prop}
 
 \begin{proof}
	 It is a consequence of Theorem \ref{Thm:eccint}, since the 
	 existence of the limit above implies $\subc = \supc = d = 
	 \frac{1} {\log2} \log\left(\lim\frac{\m_n} {\m_{2n}}\right)$.
 \end{proof}
 
 \subsection{Direct sums and tensor products of spectral triples}
 
 We study here the behaviour of noncommutative dimensions under direct
 sum and tensor product.

 \begin{Prop}\label{mu-oplus}
	Let $A$, $B$ be compact operators. Then 
	\begin{equation*}
		\m_{A\oplus B}=\m_{A}\oplus\m_{B}.
	\end{equation*}
 \end{Prop}

 \begin{proof}
	In the definition of $\m_{A}\oplus\m_{B}$, choose
	$\tilde\m_{A}$ to be the function defined on two copies of
	$\br_{+}$ which is equal to $\m_{A}$ on the first copy and to
	zero on the second.  Analogously, set $\tilde\m_{B}$ to be
	equal to $\m_{B}$ on the second copy and to zero on the first. 
	Recall that the distribution function of $\m$ is $\l_{\m}(t) := meas 
	\{x>0 : \m(x)>t\}$ (cf.  \cite{BS,FK}). We clearly have
	$\l_{\tilde\m_{A}+\tilde\m_{B}}=\l_{\m_{A}}+\l_{\m_{B}}$ and
	also $\l_{A\oplus B}=\l_{\m_{A}}+\l_{\m_{B}}$.  The thesis follows.
 \end{proof}

 \begin{Cor}\label{sumdim}
     Let $\triple_{i}=(\ca_{i},\ch_{i},D_{i})$, $i=1,2$, and
     $\triple=(\ca,\ch_{1}\oplus\ch_{2},D_{1}\oplus D_{2})$, be
     spectral triples.  Then $d(\triple)=d(\triple_{1})\vee
     d(\triple_{2})$.  The interval $[\subc(\triple), \supc(\triple)]$
     is contained in the interval
     $[\subc(\triple_{1})\wedge\subc(\triple_{2}),
     \supc(\triple_{1})\vee\supc(\triple_{2})]$.
 \end{Cor}

 \begin{proof}
	Immediately follows by Propositions \ref{estimates},
	\ref{mu-oplus} and Theorem \ref{c(oplus)}
 \end{proof}

 \medskip
 
 Let $\triple_{i}=(\ca_{i},\ch_{i},D_{i})$, $i=1,2$, be spectral
 triples.  Then their tensor product is the spectral triple
 $\triple_{1}\otimes\triple_{2}=(\ca,\ch,D)$, where $\ca :=
 \ca_{1}\otimes\ca_{2}$, $\ch:=\ch_{1}\otimes \ch_{2}$, and $D$ is
 defined in different ways according to the parity of the two triples,
 but $D^{2}$ is, up to a finite multiplicity, always equal to
 $D_{1}^{2}\otimes 1 + 1\otimes D_{2}^{2}$.

 \begin{Prop}
	With notation as above, 
	$$
	d(\triple_{1}\otimes\triple_{2}) \leq d(\triple_{1}) + 
	d(\triple_{2}).
	$$
 \end{Prop}
 \begin{proof}
	Let $\z_{D}(\a) := \sum_{n=0}^{\infty} \m_{n}(D)^{\a}$, 
	$\a\in\br$, denote the ``zeta'' function of the spectral triple 
	$(\ca,\ch,D)$, and analogously for $(\ca_{i},\ch_{i},D_{i})$, 
	$i=1,2$.  Then, if $c\in\bn$ denotes the multiplicity, and 
	$\a_{i}> d(\ca_{i},\ch_{i},D_{i})$, $i=1,2$, we have
	\begin{align*}
		\z_{D}(\a_{1}+\a_{2}) & = \sum_{n=0}^{\infty} 
		\m_{n}(D^{2})^{-(\a_{1}+\a_{2})/2} = 
		c \sum_{n=0}^{\infty} \m_{n}(D_{1}^{2}\otimes 1 + 1\otimes 
		D_{2}^{2})^{-(\a_{1}+\a_{2})/2} \\
		& = c \sum_{i,j=0}^{\infty} 
		\{\m_{i}(D_{1})^{2}+\m_{j}(D_{2})^{2}\}^{-(\a_{1}+\a_{2})/2} 
		\\
		&\leq c \sum_{i,j=0}^{\infty} 
		\m_{i}(D_{1})^{-\a_{1}}\m_{j}(D_{2})^{-\a_{2}}
		= c \z_{D_{1}}(\a_{1})\z_{D_{2}}(\a_{2}),
	\end{align*}
	which converges.  Therefore, by Theorem \ref{Thm:dimension}, we 
	get the thesis.
 \end{proof}
 
 \section{Fractals in $\br$. Classical aspects}\label{SecClassic}
 
 \subsection{Preliminaries}
 
 Let $(X,\r)$ be a metric space, and let $h:[0,\infty) \to [0,\infty)$ 
 be non-decreasing and right-continuous, with $h(0)=0$.  When 
 $E\subset X$, define, for any $\d>0$, $\ch^{h}_{\d}(E) := \inf \{ 
 \sum_{i=1}^{\infty} h(\diam A_{i}) : \cup_{i} A_{i} \supset E, \diam 
 A_{i} \leq \d \}$.  Then the {\it Hausdorff-Besicovitch (outer) 
 measure} of $E$ is defined as
 $$                                                                               
   \ch^{h} (E) := \lim_{\d\to0}\ch^{h}_{\d}(E).                                     
 $$                                                                               
 If $h(t) = t^{\a}$, $\ch^{\a}$ is called {\it Hausdorff (outer) 
 measure} of order $\a>0$.

 The number
 $$                                                                               
 d_{H}(E) := \sup \{ \a>0 : \ch^{\a}(E) = +\infty \} = \inf \{ \a>0 : 
 \ch^{\a}(E) = 0 \}
 $$                                                                               
 is called {\it Hausdorff dimension} of $E$.                                      

 Let $N_{\eps}(E)$ be the least number of closed balls of radius 
 $\eps>0$ necessary to cover $E$.  Then the numbers
 $$                                                                               
 \ubd{E} := \limsup_{\eps\to0^{+}} \frac{\log N_{\eps}(E)}{-\log 
 \eps}, \quad \lbd{E} := \liminf_{\eps\to0^{+}} \frac{\log 
 N_{\eps}(E)}{-\log \eps}
 $$                                                                               
 are called upper and lower {\it box dimensions} of $E$.                          
                                                                                    
 In case $X=\br^{N}$, setting $S_{\eps}(E) := \{ x\in\br^{N} : \r(x,E) 
 \leq \eps \}$, it is known that $\ubd{E} = N - \liminf_{\eps\to0^{+}} 
 \frac{\log \vol S_{\eps}(E)}{\log \eps}$ and $\lbd{E} = N - 
 \limsup_{\eps\to0^{+}} \frac{\log \vol S_{\eps}(E)}{\log \eps}$.  $E$ 
 is said {\it $d$-Minkowski measurable} if the following limit exists:
 $$                                                                               
 \cam_{d}(E):=\lim_{\eps\to0^{+}} \frac{\vol 
 S_{\eps}(E)}{\eps^{N-d}}\in (0,\infty).
 $$                                                                               
 This implies that the upper and lower box dimensions coincide.  The
 quantity $\cam_{d}(E)$ is called {\it $d$-Minkowski content} of $E$.
 He and Lapidus \cite{HeLa} have recently generalised that as 
 follows. If $h:[0,\infty)\to[0,\infty)$ is non-decreasing and 
 $h(0)=0$, $E$ is said {\it $h$-Minkowski measurable} if the following 
 limit exists:
 $$
 \cam_{h}(E) := \lim_{\eps\to0} \vol 
 S_{\eps}(E)\ \frac{h(\eps)}{\eps^{N}}\in (0,\infty).
 $$
 The quantity $\cam_{h}(E)$ is called {\it $h$-Minkowski content} of
 $E$.

 \subsection{Fractals in $\br$}\label{FractalsinR}
 
 By a fractal in $\br$ we mean a compact, totally disconnected subset
 of $\br$, without isolated points.  Let $F$ be such a set, and denote
 by $[a,b]$ the least closed interval containing $F$.  Then
 $[a,b]\setminus F$ is the disjoint union of open intervals
 $(a_{n},b_{n})$, which we assume ordered in such a way that
 $\{b_{n}-a_{n}\}_{n\in\bn}$ is a decreasing sequence.  Notice that
 $F$ is determined by the sequence of intervals
 $\{(a_{n},b_{n})\}_{n\in\bn}$.  Then $F$ has Lebesgue measure zero
 $iff$ $\sum_{n=1}^{\infty} (b_{n}-a_{n}) = b-a$, and, in that case,
 (cf.  e.g. \cite{Tricot})
 \begin{equation}\label{dimlin}
	 \ubd{F}=\limsup_{n\to\infty} \frac{\log n}{|\log(b_{n}-a_{n})|}.
 \end{equation}

 We will be interested in fractals constructed out of a family 
 $\{w_{ni}:i=1,\ldots p_{n}, n\in\bn\}$ of contracting similarities 
 of $\br$, with dilation parameters $\l_{ni}$, such that 
 \begin{itemize}
 	\itm{i} $w_{ni}([a,b]) \subset [a,b]$
	\itm{ii} $w_{ni}([a,b]) \cap w_{nj}([a,b]) = \emptyset, i\neq j,\ 
	n\in\bn$
	\itm{iii} $\bigcup_{i=1}^{p_{n}} w_{ni}(\{a,b\}) \supset \{a,b\}$, 
	$n\in\bn$.
 \end{itemize}
 For any $n\in\bn$, set $w_{n}(\O) := \bigcup_{i=1}^{p_{n}} w_{ni}(\O)$, 
 $\O\subset\br$, and $W_{n}:= w_{1}\circ 
 w_{2}\circ\cdots\circ w_{n}$. Then $\{W_{n}([a,b])\}$ is a 
 decreasing sequence of compact sets, containing $\{a,b\}$. Denote by 
 $F$ its intersection. Then
 
 \begin{Prop}
 	$F$ is a fractal in $\br$. It has Lebesgue measure zero iff 
 	$$\prod_{n}\left(\sum_{i=1}^{p_{n}}\l_{ni}\right) = 0.$$
 \end{Prop}
 
 We call the fractals described above {\it limit fractals} (cf.  
 \cite{GuIs8,GuIs10} for alternate, more general definitions).  If 
 $p_{n}=p$, for all $n\in\bn$, and the similarity parameters $\l_{ni}$ 
 do not depend on $n$, $F$ is a self-similar fractal \cite{Hutch}.  If the 
 similarity parameters $\l_{ni}$ do not depend on $i$, $F$ is called a 
 {\it translation fractal} (cf.  \cite{KaSa}). Observe that for a 
 translation fractal the condition $(ii)$ above implies $p_{n}\l_{n}<1$
 
 The fractal is called {\it symmetric}, if $w_{n}([a,b]) =
 \bigcup_{i=1}^{p_{n}} [a+(i-1)d_{n},a+(i-1)d_{n}+\l_{n}]$, where
 $d_{n} := \frac{(b-a)-p_{n}\l_{n}}{p_{n}-1}$.  In this case $F$ is
 uniquely determined by the sequences $\{p_{n}\}$, $\{\l_{n}\}$.
 
 \subsection {Symmetries of limit fractals} 
 
 Let us denote by $F_{n}$ the set $\bigcap_{k=0}^{\infty} w_{n+1}\circ
 w_{n+2}\circ\cdots\circ w_{n+k}([a,b])$.  We clearly have $F =
 w_{1}\circ w_{2}\circ\cdots\circ w_{n}(F_{n})$.  Therefore, if $\s$
 denotes a multiindex of length $|\s|=n$, and $w_{\s}:=
 w_{1\s(1)}\circ w_{2\s(2)}\circ\cdots\circ w_{n\s(n)}$, we have $F =
 \bigcup_{|\s|=n}w_{\s}(F_{n})$, with disjoint union.
 
 We call the similarity maps $w_{\s'}\circ 
 w_{\s}^{-1}:w_{\s}(F_{n})\mapsto w_{\s'}(F_{n})$, $|\s|=|\s'|=n$, 
 $n\in\bn$, {\it generating symmetries} of the limit fractal $F$.
 
 Observe that if the fractal is a translation fractal the generating 
 symmetries are indeed isometries.
 
 Let us consider a triple $(\O_{1},\O_{2},S)$ where $\O_{1},\O_{2}$
 are (relatively) open subsets of $F$ and $S$ is a one-to-one similarity
 with scaling parameter $\l$ between $\O_{1}$ and $\O_{2}$.  We say
 that a measure $\m$ on $F$ is homogeneous of order $\a>0$ for the
 triple $(\O_{1},\O_{2},S)$ if $\m(\O_{2}) = \l^{\a}\m(\O_{1})$.
 
 \begin{Prop}\label{homogen-meas}
	Let $F$ be a limit fractal.  Then, for any $\a\in(0,1)$, there
	is a unique probability measure $\m_{\a}$, with support $F$,
	homogeneous of order $\a$ w.r.t. the generating symmetries of
	the fractal.  All these measures are distinct, unless $F$ is a
	translation fractal, in which case they all coincide.
 \end{Prop}
 
 \begin{proof}
 	For any $n$, the homogeneity condition uniquely determines the 
 	measure of the sets $w_{\s}(F_{n})$, $|\s|=n$. Indeed, if $w_{\s}$ 
 	has similarity parameter $\l_{\s}$, 
	\begin{align*}
		1&=\m_{\a}(F)=\sum_{|\s'|=n}\m_{\a}(w_{\s'}(F_{n}))\\
		&=\sum_{|\s'|=n}\m_{\a}(w_{\s'}\circ w_{\s}^{-1}(w_{\s}(F_{n})))\\
		&=\sum_{|\s'|=n}(\l_{\s'}\l_{\s}^{-1})^{\a}\m_{\a}(w_{\s}(F_{n}))\\
		&=\l_{\s}^{-\a}\m_{\a}(w_{\s}(F_{n}))\sum_{|\s'|=n}(\l_{\s'})^{\a}
	\end{align*}
	namely 
	$$
	\m_{\a}(w_{\s}(F_{n})) 
	=\l_{\s}^{\a}\left(\sum_{|\s'|=n}(\l_{\s'})^{\a}\right)^{-1}.
	$$
	The measure uniquely extends to the sigma-algebra generated by 
	these sets, which clearly coincides with the family of Borel 
	subsets of $F$.  The second statement is obvious.
 \end{proof}
  
 \begin{rem}\label{TranslationFractals}
     It has been proved in \cite{KaSa} that, when $F$ is a translation
     fractal in $\br$, there is a gauge function $h$ such that the
     corresponding Hausdorff-Besicovitch measure $\ch^{h}$ is
     non-trivial on $F$.  Since any Hausdorff-Besicovitch measure is
     isometry invariant, it satisfies the hypotheses of the previous
     proposition, hence $\ch^{h}|_F$ coincides (up to a constant) with
     the homogeneous measure $\m$.
 \end{rem}
 
 \section{Fractals in $\br$. Noncommutative 
 aspects.}\label{NonCommFractals}

 \subsection{The lacunary spectral triple}

 Let $F$ be a fractal in $\br$, namely a compact, totally disconnected
 subset of $\br$, without isolated points.  Now we introduce a
 ``lacunary'' spectral triple for the fractal $F$, namely a spectral
 triple completely determined by the ``lacunae'' of $F$, hence in
 particular canonically associated to $F$.  Amendments to this
 spectral triple will be discussed below.  Let $a,b,a_n,b_n$ be as in
 subsection \ref{FractalsinR}, and denote by $I_{n}$ the lacuna
 $(a_{n},b_{n})$.  Set $\ch_{\ell} = \oplus_{n=1}^{\infty}
 \ch(I_{n})$, $D_{\ell}= \oplus_{n=1}^{\infty} D(I_{n})$, where
 \begin{align}
	 \ch(I)&:= \ell^{2}(\partial I), \label{HilbertInterval}\\
	 D(I) &:= \frac{1}{|I|} 
	 \begin{pmatrix}
		 0 & 1\\
		 1 & 0
	 \end{pmatrix}.\label{DiracInterval}
 \end{align}
  Consider the action of $C(F)$ on 
 $\ch_{\ell}$ by left multiplication: $(f\xi)(x)=f(x)\xi(x)$, $x\in 
 \cd_{\ell}:=\{a_n,b_n : n\in\bn\}$, and define $\ca := Lip(F)$.  Then
 
 \begin{Thm}\label{linfrac}
 	\itm{i} $(\ca,\ch_{\ell},D_{\ell})$ is a spectral triple
	\itm{ii} the characteristic values of $D_{\ell}^{-1}$ are the numbers 
	$b_{n}-a_{n}$, $n\in\bn$, each with multiplicity 2.
 
 	If $F$ is Minkowski measurable, and has box dimension $d\in(0,1)$, 
 	then
	\itm{iii} $|D_{\ell}|^{-d}\in \cl^{1,\infty}$
	\itm{iv} $\Tr_{\o}(|D_{\ell}|^{-d}) = 2^{d}(1-d) \cam_{d}(F)$.
 \end{Thm}
 \begin{proof}
 	It is due to Connes \cite{Co}, using results of Lapidus and Pomerance, 
 \cite{LaPo}.
 \end{proof}
  
 Making use of recent results of He and Lapidus \cite{HeLa}, we can 
 improve on the previous Theorem. Recall from \cite{HeLa} that the 
 family of gauge functions $G_{d}$, for $d\in(0,1)$, consists of the 
 functions $h:(0,\infty)\to(0,\infty)$ which are continuous, strictly 
 increasing, with $\lim_{x\to0} h(x) = 0$, $\lim_{x\to\infty} h(x) = 
 \infty$, and satisfy
 $$
 \lim_{x\to0}\frac{h(tx)}{h(x)} = t^{d}
 $$
 uniformly in $t$ on any compact subset of $(0,\infty)$, and one more 
 condition (H3), which won't be needed in the following. Then, 
 setting $g(x) := h^{-1}(1/x)$, $x>0$, we have
 
 \begin{Thm}\label{h-Minkowski}
     Let $d\in(0,1)$, and $h\in G_{d}$, and assume $F$ is 
     $h$-Minkowski measurable. Then
     \itm{i}  the function $g^{d}$ is eccentric, so it gives rise to 
     a singular trace 
     $$
     \t_{h,\o}(a) = \Lim_{\omega} \left(\frac{S_{a}(n)} { S_{g^{d}}(n) }\right)
     $$
     \itm{ii} $d=d(\ca,\ch_{\ell},D_{\ell})$ and is the unique exponent of 
     singular traceability of $D_{\ell}^{-1}$
     \itm{iii} $\t_{h,\o}(|D_{\ell}|^{-d}) = 2^{d}(1-d) \cam_{h}(F)$, and is 
     therefore independent of the state $\o$.
 \end{Thm}
 \begin{proof}
     $(i)$ Recall from \cite{HeLa}, Theorems 2.4 and 2.5, that $F$ is
     $h$-Minkowski measurable $iff$ there is $L>0$ such that
     $b_{n}-a_{n} \sim Lg(n)$, $n\to\infty$, and in this case 
     $\cam_{h}(F) = \frac{2^{1-d}L^{d}}{1-d}$.
     
	 Besides, it follows from \cite{HeLa}, Lemma 3.1 that
     \begin{equation*}
	 \lim_{z\to\infty}\frac{g(tz)}{g(z)} = \frac{1}{t^{1/d}},
     \end{equation*}
     for any $t>0$. Therefore 
     \begin{equation*}
	 \lim_{z\to\infty} \frac{\log g(z)}{\log 1/z} = \frac{1}{d}
     \end{equation*}
     which shows that $g^{d}$ is eccentric.
     
     $(ii)$ First observe that
     \begin{equation}\label{eq:modified.reg}
	 \lim_{z\to\infty}\frac{g(tz+a)}{g(z)} = \frac{1}{t^{1/d}},
     \end{equation}
     for any $t>0$, $a\in\br$. Indeed, for any $\eps\in(0,t)$, there is 
     $z_{\eps}>0$ such that $(t-\eps)z \leq tz+a \leq (t+\eps)z$, 
     $z>z_{\eps}$, so that
     \begin{align*}
	 \frac{1}{(t+\eps)^{1/d}} = 
	 \lim_{z\to\infty}\frac{g((t+\eps)z)}{g(z)} \leq 
	 \lim_{z\to\infty}\frac{g(tz+a)}{g(z)} \leq 
	 \lim_{z\to\infty}\frac{g((t-\eps)z)}{g(z)} = \frac{1}{(t-\eps)^{1/d}},
     \end{align*}
     and the thesis follows from the arbitrariness of $\eps$.
         
     Let us now denote by $\m_{n}$ the $n$-th characteristic value of $D_{\ell}^{-1}$. 
     Because of the previous Theorem, $\m_{2n-1} = \m_{2n} = 
     b_{n}-a_{n}$, so that
     \begin{align*}
	 \lim_{n\to\infty} \frac{\m_{2n}}{\m_{4n}} & =
	 \lim_{n\to\infty} \frac{g(n)}{g(2n)} = 2^{1/d} \\
	 \lim_{n\to\infty} \frac{\m_{2n-1}}{\m_{4n-2}} & =
	 \lim_{n\to\infty} \frac{g(n)}{g(2n-1)} = 2^{1/d},
     \end{align*}
	 where the last equality follows from (\ref{eq:modified.reg}).  
	 Therefore $\lim_{n\to\infty} \frac{\m_{n}}{\m_{2n}} = 2^{1/d}$, 
	 and, by Proposition \ref{Prop:unique} and its proof, we conclude.
     
     $(iii)$ Let us first observe that
     $$
     \exists\ \lim_{t\to\infty} \frac{\m(t)}{g(t)} = \a\in 
     [0,\infty]\ \iff\ \exists\ \lim_{n\to\infty} \frac{\m(2n)}{g(2n)} = \a\in 
     [0,\infty].
     $$
     Indeed, for any $t>0$, there is $n\in\bn$ such that 
     $t\in(2n-2,2n]$, so that
     $$
     \frac{g(2n)}{g(2n-2)}\ \frac{\m(2n)}{g(2n)} = 
     \frac{\m(2n)}{g(2n-2)} \leq \frac{\m(t)}{g(t)} \leq 
     \frac{\m(2n-2)}{g(2n)} = \frac{\m(2n-2)}{g(2n-2)}\ \frac{g(2n-2)}{g(2n)}
     $$
     and the thesis follows from (\ref{eq:modified.reg}).
     
     Now assume that $\int_{0}^{\infty} g(t)^{d}\ dt =\infty$. Then
     \begin{align*}
	 \lim_{t\to\infty} \frac{\int_{0}^{t} \m(s)^{d}\
	 ds}{\int_{0}^{t} g(s)^{d}\ ds} & = \left(\lim_{t\to\infty}
	 \frac{\m(t)}{g(t)}\right)^{d} = \left(\lim_{n\to\infty}
	 \frac{\m(2n)}{g(2n)}\right)^{d} \\
	  &= \left(\lim_{n\to\infty} \frac{b_{n}-a_{n}}{g(n)}\ 
	  \frac{g(n)}{g(2n)}\right)^{d} 2L^{d} = 2^{d}(1-d)\cam_{h}(F).
     \end{align*}
     We can proceed in a similar way if $\int_{0}^{\infty} g(t)^{d}dt <\infty$. 
 \end{proof}    
 
  Even if $F$ is not $h$-Minkowski measurable, we have that, by Theorem 
 \ref{Thm:trace}, any singular traceability exponent gives rise to a 
 trace state on the C$^{*}$-algebra of continuous functions on the 
 fractal, namely to a probability measure on the fractal.  In 
 particular,
 
 \begin{Thm}\label{d=box}
	 $(i)$ For any singular traceability exponent $s$ for
	 $|D_{\ell}|^{-1}$ we get a Hausdorff-Be\-si\-co\-vitch functional on
	 the spectral triple, giving rise to a probability measure $\mu$
	 on $F$.\\
	 $(ii)$ Let $F$ have zero Lebesgue measure.  Then
	 $d(\ca,\ch_{\ell},D_{\ell}) = \ubd{F}$.  Therefore, if
	 $\ubd{F}\ne0$, we get a corresponding measure on $F$.
 \end{Thm}
 \begin{proof}
	 $(i)$ follows from Theorem \ref{Thm:trace} (cf. Definition 
	 \ref{Dfn:dimensions}) and Riesz Theorem.
	 \\
	 $(ii)$ follows by equation (\ref{dimlin}) and 
	 Theorem \ref{linfrac} $(i),(ii)$.  
 \end{proof}

 \subsection{The reconstruction of the metric}
 
 First we discuss the ``lacunary'' metric on the fractal, namely the 
 metric on $F$ determined {\it \`a la Connes} via the lacunary 
 spectral triple. As explained below, such metric does not coincide 
 in general with the original one.
 
 Let us first compute $\|[D_{\ell},f]\|$. Observe that, setting 
 $I=(x,y)$,
 $$
 \|[D(I),f|_{\partial I}]\|=\frac{1}{|I|}
 \left\|\left[\left(\begin{matrix}0&1\\1&0\end{matrix}\right),
 \left(\begin{matrix}f(x)&0\\0&f(y)\end{matrix}\right)\right]\right\|
 =\left|\frac{f(y)-f(x)}{y-x}\right|.
 $$
 Therefore we have
 \begin{equation}\label{metricIneq}
	 \|[D_{\ell},f]\|=\sup_{n}\|[D(I_{n}),f|_{\partial I_{n}}]\|=
	 \sup_{n}\left|\frac{f(b_{n})-f(a_{n})}{b_{n}-a_{n}}\right|
	 \leq\|f\|_{Lip(F)}.
 \end{equation}
 
 \begin{Thm}\label{lacunary metric}
	 Let $F$ be a compact, totally disconnected subset of $\br$ with no 
	 isolated points. Then the lacunary metric
	 \begin{equation}\label{ConnesMetric}
		 d_{\ell}(x,y)=\sup\{|f(y)-f(x)|:f\in C(F), \|[D_{\ell},f]\|\leq1\}
	 \end{equation}
	 coincides with the one induced by the metric on $\br$ if and only 
	 if $F$ has Lebesgue measure zero.
 \end{Thm}
 
 \begin{proof}
	Assume $F$ has Lebesgue measure zero. Now, for any $f\in C(F)$, we 
	denote by $\tilde{f}$ the continuous function on $[a,b]$ coinciding 
	with $f$ on $F$ and linear on any interval $[a_{n},b_{n}]$. For any 
	pair of points $x<y$ in $F$, let us denote by $\cai(x,y)$ the family 
	of lacunary intervals $I_{n}$ which are subsets of $[x,y]$. Then
	\begin{align*}
		|f(y)-f(x)|&
		=\left|\sum_{I\in\cai(x,y)}\int_{I}\tilde{f}'(t)dt\right|\\
		&
		\leq \sum_{I\in\cai(x,y)}|I|\|[D(I),f|_{\partial I}]\|
		\leq |y-x| \|[D_{\ell},f]\|.
	\end{align*}
	Comparing the previous inequality with (\ref{metricIneq}) we get 
	$\|f\|_{Lip(F)}=\|[D_{\ell},f]\|$, namely the equality $d_{\ell}=d$.
	\\
	Conversely, assuming $F$ has positive Lebesgue measure, let $f$ 
	be the restriction to $F$ of the primitive of the characteristic
	function of $F$. Clearly $\|[D_{\ell},\lambda f]\|=0$ for any 
	$\lambda>0$, hence $d_{\ell}(x,y)=+\infty$ for any pair $x,y$ in $F$ 
	which are not boundary of the same lacuna.
 \end{proof}
 
 Connes proposed us an emendation of the lacunary spectral triple in
 order to reproduce the original distance also in the case of positive
 Lebesgue measure.  In the case of the Cantor middle third set, the idea is to
 add to the lacunary intervals also the images of the interval $[0,1]$
 under the similarity maps \cite{Cpc}. 
 
 For a general compact, totally disconnected fractal $F$, the idea of Connes
 may be generalized as follows:
 
 Assign the family $\cf_n$ of (closed) filled intervals of level $n$ and the 
 family $\cl_n$ of (open) lacunary intervals of level $n$ in such a way that
 \begin{itemize}
 	\item $F=\cap_n\cup_{I\in\cf_n}I$,
	\item $\cl_n\subseteq\cl_{n+1}$,
	\item $\cai_n:=\cf_n\cup\cl_n$ form a finite partition of 
	$[a,b]$ for any $n$.
 \end{itemize}
 
 Then, setting $\cai=\cup_n\cai_n$, and, according to the notation 
 in subsections 3.2 and 3.3, 
 $\ch=\oplus_{I\in\cai}\ch(I)$, $D=\oplus_{I\in\cai}D(I)$, we get that 
 $(\ca,\ch,D)$ is a spectral triple, $\ca$ being the $^{*}$-algebra 
 of Lipschitz functions. Moreover,
 
 \begin{Thm}\label{NewTriple}
	 Let $F$ be a compact, totally disconnected subset of $\br$ with
	 no isolated points. Then the spectral triple $(\ca,\ch,D)$ 
	 reconstructs the original distance on $F$.
 \end{Thm}
 
 \begin{proof}
	For any pair $x<y$ in the boundary of some lacunae, there exists 
	a $k$ such that both $x$ and $y$ belong to the boundary of some 
	interval of level $k$. Therefore, setting $\cai_{k}(x,y)$ for the family 
	of intervals of level $k$ which are subsets of $[x,y]$, we have
	\begin{align*}
		|f(y)-f(x)|&
		\leq \sum_{I\in\cai_{k}(x,y)}|I|\|[D(I),f|_{\partial I}]\|
		\leq |y-x| \|[D,f]\|.
	\end{align*}
	Since $F$ is totally disconnected, $x$ and $y$ vary in a dense
	subset of $F$, therefore, by continuity, the previous
	inequality holds for any pair $x,y\in F$, giving
	$\|f\|_{Lip(F)}\leq\|[D,f]\|$.  On the other hand, as in
	(\ref{metricIneq}), the converse inequality holds too, hence
	the result follows.
 \end{proof}

 \subsection{A spectral triple for limit fractals}

 The spectral triple described in the previous subsection depends on
 the choice of the filled and lacunary intervals of level $n$.  Of
 course, one may either select the filled intervals of level $n$
 first, and then the lacunae as the connected components of the
 complement, or the converse.
 The first choice appears very natural in the case of limit fractals,
 therefore we shall adopt this point of view, limiting our further
 analysis to this family.
 
 \begin{Dfn}\label{Dfn:newtriple}
     Let $F$ be a limit fractal, with similarities $w_{n,i}$.
     We set $\cf_n$ to be $w_{\s}[a,b]$, where $\s$ varies in the set of
     multi-indices of length $n$.
 
     The Dirac operator is a direct sum of
     the lacunary Dirac $D_{\ell}$ and the Dirac 
     $D_{f}=\oplus_{n\in\bn,I\in\cf_n}D(I)$ 
     acting on the Hilbert space $\ch_{f}=\oplus_{n\in\bn,I\in\cf_n}\ch(I)$.
 
     We choose $\ca$ to be the $^{*}$-algebra of Lipschitz functions,
     acting on $\ch$ by pointwise multiplication.
 \end{Dfn}

 \begin{rem}
     \item[(1)] The spectral triple $(\ca,\ch,D)$ reconstructs the
     original distance on $F$, by Theorem \ref{NewTriple}. Any 
     singular traceability exponent, namely any number between 
     $\subc(\ca,\ch,D)$ and $\supc(\ca,\ch,D)$ gives rise to a 
     Hausdorff-Besicovitch functional on $\ca$, by Theorems 
     \ref{Thm:SingTracExp} and \ref{Thm:trace}.
     
     \item[(2)] $(\ca,\ch_{f},D_{f})$ is a spectral triple, indeed
     it is exactly the spectral triple we proposed in \cite{GuIs8} for
     limit fractals in $\br^{n}$.
     
     \item[(3)] Concerning the second choice, namely defining the
     intervals $\cai$ selecting the lacunae first, one may e.g. call
     $\lambda_n$ the values of the lengths of the lacunae arranged in
     decreasing order, and then $\cl_n$ the lacunae of length lower
     equal than $\lambda_n$.  While this choice is completely
     canonical, even the analysis of self-similar fractals is far less
     obvious than the corresponding one for the lacunary spectral
     triple.  However, the spectral triple corresponding to such a
     choice coincides with the one in Definition \ref{Dfn:newtriple}
     in the case of uniformly generated symmetric fractals.
 \end{rem}

 
 In the following we shall prove some results on the dimensions and
 measures associated with the spectral triple of a limit fractal.  The
 results are generally stated for the triple $(\ca,\ch,D)$, but hold
 also for the ``lacunary'' and ``filled'' spectral triples.  Indeed we
 shall prove such properties for the latter triples, then showing that
 they remain valid for the direct sum of the Dirac operators.
 
 \begin{Thm}\label{Meas*LimFrac}
	Let $F$ be a limit fractal, $s$ a singular traceability exponent
	for $|D|^{-1}$, $\tau_{\omega}$ a corresponding singular trace. 
	Then, for any continuous function $f$ on $F$,
	\begin{equation}\label{eqmeas}
		\tau_{\omega}(f|D|^{-s})=\int_{F}f\ d\mu_{s},
	\end{equation}
	where $\mu_{s}$ is the measure introduced in Proposition
	\ref{homogen-meas}.  In particular continuous functions are
	measurable, namely the integral of continuous functions does not
	depend on the generalized limit procedure $\omega$.
 \end{Thm}

 Before proving the theorem, we give a corollary which follows 
 immediately from Remark \ref{TranslationFractals}
 
 \begin{Cor}
     If $F$ is a translation fractal, all the Hausdorff-Besicovitch
     functionals on $\ca$ described above give rise to the same
     measure, which is indeed the restriction to $F$ of the
     Hausdorff-Besicovitch measure constructed in \cite{KaSa}.
 \end{Cor}
 
 Let us first discuss the statement above for the lacunary Dirac 
 operator.
 
 \begin{Lemma}
	 The probability measure $\mu$ on $F$ associated with a singular
	 traceability exponent $s$ for $|D_{\ell}|^{-1}$ has the following
	 property:
	 \begin{equation}\label{homogeneity}
		 \m(\O_{2})=\l^{s}\m(\O_{1})
	 \end{equation}
	 where $\O_{1}$, $\O_{2}$ are relatively open subsets of $F$ 
	 related by a similarity of parameter $\l$.
 \end{Lemma}
 \begin{proof}
	 Observe that given any $(\O_{1},\O_{2},S)$, where 
	 $\O_{1},\O_{2}$ are clopen sets in $F$ and $S:\O_{1}\to \O_{2}$ 
	 is a one-to-one similarity of parameter $\l$, we have that the 
	 sizes of the lacunae in $\O_{1}$ are multiples of the sizes of the 
	 lacunae in $\O_{2}$ with scaling $\l$.  Therefore, by Theorem 
	 \ref{linfrac} $(ii)$, the eigenvalues of $\chi_{\O_{1}}|D_{\ell}|^{-s}$ 
	 are multiples of the eigenvalues of $\chi_{\O_{2}}|D_{\ell}|^{-s}$ with 
	 scaling $\l^{s}$.  As a consequence the measure of $\O_{1}$ is 
	 equal to $\l^{s}$ times the measure of $\O_{2}$.  This clearly 
	 extends to pairs of open sets $\O_{1},\O_{2}$.  \\
 \end{proof}
 
 Previous lemma can be void for general fractals, namely there may be
 no non-trivial triples $(\O_{1},\O_{2},S)$.  However, if limit
 fractals are concerned, the generating symmetries determine the
 measure $\mu_{s}$ introduced in Proposition \ref{homogen-meas}, 
 hence we have proved 
 
 \begin{Prop}
	 Equation (\ref{eqmeas}) holds for the lacunary Dirac.
 \end{Prop}
 
 The proof for the filled Dirac operator is analogous.  Indeed the
 following lemma is a direct consequence of the definition of $D_{f}$.
 
 \begin{Lemma}
	 The eigenvalues of $D_{f}$ are the numbers $(b-a)\lambda_{\s}$,
	 each with multiplicity two, where we have set $\lambda_{\s} =
	 \prod_{i=1,\dots,|\s|}\lambda_{i,\s_{i}}$.
 \end{Lemma}
 
 Then, if $\O_{1},\O_{2}$ are clopen sets in $F$ related by a
 generating symmetry of $F$, the eigenvalues of $\chi_{\O_{1}}
 |D_{f}|^{-s}$ are multiples of the eigenvalues of $\chi_{\O_{2}}
 |D_{f}|^{-s}$ with scaling $\l^{s}$. Reasoning as above, we have

 \begin{Prop}
	 Equation (\ref{eqmeas}) holds for the filled Dirac.
 \end{Prop}
 
 \begin{proof} {\it (of Theorem \ref{Meas*LimFrac}).} 
	 If $\O_{1},\O_{2}$ are clopen sets in $F$ related by a generating
	 symmetry of $F$, the eigenvalues of $\chi_{\O_{1}} |D|^{-s}$ are
	 multiples of the eigenvalues of $\chi_{\O_{2}} |D|^{-s}$ with
	 scaling $\l^{s}$, since this property holds for the two direct
	 summands $D_{\ell}$ and $D_{f}$. The result then follows.
 \end{proof}

 \begin{rem}
	 Let us note that in spite of the fact that equation (\ref{eqmeas})
	 holds for the three Dirac operators, the singularity exponents for 
	 the different Dirac's are different in general. They will coincide 
	 however for uniformly generated symmetric fractals.
 \end{rem}
 
 
 \begin{Thm}\label{Thm:self-similar}
	Let $F$ be a self-similar fractal, and $d\in(0,1)$ its
	Hausdorff dimension.  Then $d$ is the unique exponent of
	singular traceability for $D^{-1}$, and the Hausdorff
	functional on the spectral triple gives rise to the
	$d$-dimensional Hausdorff measure on $F$, up to a
	multiplicative constant.  In particular the commutative and
	noncommutative Hausdorff dimensions coincide.
 \end{Thm}
 \begin{proof}
     Set $\cd_{\ell}:=\{a_n,b_n : n\in\bn\}$, $\cd_{f}:=
     \bigsqcup_{\s\in\S^{*}} \{\s(a),\s(b)\}$, where $\bigsqcup$
     denotes disjoint union, and $\S^{*}$ is the set of all
     multi-indices.  Define the following operators on
     $\ell^{2}(\cd_{\ell})$:
     $$
     S_{\ell,j}\xi (b) :=
     \begin{cases}
	 \xi(w_{j}^{-1}(b)) & b\in w_{j}\cd_{\ell} \\
	 0 & b\not\in w_{j}\cd_{\ell},
     \end{cases}
     $$
     and analogously for $S_{f,j}$ on $\ell^{2}(\cd_{f})$,
     $j=1,\ldots,p$.  Then $S_{\ell,j}$ and $S_{f,j}$ are isometries
     and $|D_{\ell}|^{-s} = \sum_{j=1}^{p} \l_{j}^{s} S_{\ell,j}
     |D_{\ell}|^{-s} S_{\ell,j}^{*}$, and an analogous formula for
     $|D_{f}|^{-s}$.  Hence, with $S_{j} := S_{\ell,j} \oplus
     S_{f,j}$, we obtain $|D|^{-s} = \sum_{j=1}^{p} \l_{j}^{s} S_{j}
     |D|^{-s} S_{j}^{*}$.  Therefore, if $s$ is an exponent of
     singular traceability for $|D|^{-1}$, the corresponding
     Hausdorff-Besicovitch functional is homogeneous of order $s$. 
     This implies that $s$ coincides with $d$, namely $d$ is the
     unique exponent of singular traceability.  We now prove that the
     $d$-dimensional Hausdorff functional corresponds to the
     $d$-dimensional Hausdorff measure. Let us compute the zeta 
     functions of $D_{\ell}$ and $D_{f}$ separately:
     \begin{align*}
	 \z_{f}(s) & := Tr(|D_{f}|^{-s}) = 2(b-a)\sum_{\s}\l_{\s}^{s} 
	  = 2(b-a) \sum_{n=0}^{\infty} \sum_{|\s|=n} 
	 \prod_{j=1}^{n}\l_{\s(j)}^{s} \\
	 & = 2(b-a) \sum_{n=0}^{\infty} \left( \sum_{j=1}^{p} \l_{j}^{s} 
	 \right)^{n} 
	  = \frac{2(b-a)}{1- \sum_{j=1}^{p} \l_{j}^{s}},
     \end{align*}
     so that
     $$
     \lim_{s\to d} (s-d)\z_{f}(s) = \frac{2(b-a)}
     {\sum_{j=1}^{p} \l_{j}^{d}\log(1/ \l_{j})},
     $$
     which, using \cite{Co}, Proposition IV.2.$\beta$.4, implies that
     the Hausdorff functional is non-trivial.  As for $D_{\ell}$,
     denoting by $c_{1},\ldots,c_{p-1}$ the lengths of the connected
     components of $[a,b]\setminus \cup_{j=1}^{p} w_{j}([a,b])$, we
     obtain
     \begin{equation*}
	 \z_{\ell}(s)  := Tr(|D_{\ell}|^{-s}) =
	 2\sum_{j=1}^{p-1}c_{j}^{s}\sum_{\s}\l_{\s}^{s} 
	  = \frac{2\sum_{j=1}^{p-1}c_{j}^{s} } {1- \sum_{j=1}^{p}
	 \l_{j}^{s} } ,
     \end{equation*}
     so that
     $$
     \lim_{s\to d} (s-d)\z_{\ell}(s) =
     \frac{2\sum_{j=1}^{p-1}c_{j}^{d}}
     {\sum_{j=1}^{p} \l_{j}^{d}\log(1/ \l_{j})},
     $$
     which, using \cite{Co}, Proposition IV.2.$\beta$.4, implies that
     the Hausdorff functional is non-trivial. The thesis follows from 
     the fact that $Tr(|D|^{-s}) = \z_{f}(s) + \z_{\ell}(s)$, and 
     Theorem \ref{Meas*LimFrac}.
  \end{proof}
 
 Also in the case of symmetric fractals we have a formula for the
 noncommutative Hausdorff dimension.  We recall that symmetric
 fractals (with convex hull $[0,1]$) are determined by two sequences
 $\{p_{n}\}$, $\{\l_{n}\}$, where $p_{n}$ is a natural number greater
 or equal than 2 and $p_{n}\l_{n}<1$.  We say that the fractal is
 uniformly generated if $\sup_{n}p_{n}<\infty$ and
 $\sup_{n}p_{n}\l_{n}<1$.

 \begin{Thm}
	Let $(\ca,\ch,D)$ be the spectral triple associated with a
	uniformly generated symmetric fractal $F$, where the
	similarities $w_{n,i}$, $i=1,\dots,p_{n}$ have scaling
	parameter $\l_{n}$.  Then
	$$
	d(\ca,\ch,D)
	=\limsup_{n}\frac{\sum_1^n\log p_k}{\sum_1^n\log 1/\l_k},
	$$
	and such dimension coincides with the upper box dimension of $F$.
 \end{Thm}
 \begin{proof}
     The thesis will follow from the next two Propositions, Corollary 
     \ref{sumdim} and Theorem \ref{d=box}, $(ii)$.
 \end{proof}
  
 \begin{Prop}
	Let $F$ be a uniformly generated symmetric fractal as before.
	Then
	$$
	d(\ca,\ch_{\ell},D_{\ell})
	=\limsup_{n}\frac{\sum_1^n\log p_k}{\sum_1^n\log 1/\l_k}.
	$$
 \end{Prop}
 \begin{proof} 
     It is not restrictive to assume $a=0$, $b=1$.  Then the
     eigenvalues of $|D_{\ell}|^{-1}$ are given by
     \begin{equation}\label{eigenvalues}
	 \tilde\La_{k}=\frac{1-p_{k+1}\l_{k+1}}{p_{k+1}-1}\prod_{j=1}^{k}\l_{j}
     \end{equation}
     with multiplicity $2\tilde P_{k} = 2(p_{k+1}-1)
     \prod_{j=1}^{k}p_{j}$, $k\in\bn\cup\{0\}$.  Therefore,
	 $$
	 Tr(|D_{\ell}|^{-\a}) = 2\sum_{k=0}^{\infty} 
	 (1-p_{k+1}\l_{k+1})^{\a}(p_{k+1}-1)^{1-\a}\prod_{i=1}^{k} 
	 p_{i}(\l_{i})^{\a}.
	 $$
	 Setting $\La_{n} := \prod_{k=1}^{n} \l_{k}$, $P_{n} := 
	 \prod_{k=1}^{n} p_{k}$, we obtain
	 \begin{align*}
		Tr & (|D_{\ell}|^{-\a})= 2\sum_{k=0}^{\infty} 
		(1-p_{k+1}\l_{k+1})^{\a}(p_{k+1}-1)^{1-\a}P_{k}\La_{k}^{\a} \\
		& = 2\sum_{k=0}^{\infty} 
		(1-p_{k+1}\l_{k+1})^{\a}(p_{k+1}-1)^{1-\a} \exp\left(\log 
		P_{k}\left(1-\a\frac{\log 1/\La_{k}}{\log P_{k}}\right)\right).
	 \end{align*}
	 
	 By the $n$-th root criterion for series, the series 
	  diverges/converges if 
	 $$
		\limsup_{k}\left(\frac {(1-p_{k+1}\l_{k+1})^{\a}}{(p_{k+1}-1)^{\a-1}} 
		\exp\left(\log P_{k}\left(1-\a\frac{\log 1/\La_{k}} 
		{\log P_{k}}\right)\right)\right)^{1/k} 
		\gtrless 1,
	 $$
	 namely, since by the uniform generation assumption 
	 $\lim_{k}\left(1-p_{k}\l_{k}\right)^{1/k}=1$ and 
	 $\lim_{k}\left(p_{k}-1\right)^{1/k}=1$, if
	 $$
		\limsup_{k} 
		\frac {\log P_{k}}{k}\left(1-\a\frac{\log1/\La_{k}}{\log P_{k}}\right) 
		\gtrless 0,
	 $$
	 and, finally, if 
	 $$
		\limsup_{k}\frac{\log P_{k}}{\log 1/\La_{k}}
		\gtrless \a,
	 $$
	which implies that 
	$\limsup_{k}\frac{\log P_{k}}{\log 1/\La_{k}}$ is the abscissa 
	of convergence of the zeta function of $|D_{\ell}|^{-1}$, hence the spectral 
	dimension by Theorem \ref{unisingtrac} ($i$).
 \end{proof}

 \begin{Prop}
	Let $F$ be a uniformly generated symmetric fractal as before.
	Then
	$$
	d(\ca,\ch_{f},D_{f})
	=\limsup_{n}\frac{\sum_1^n\log p_k}{\sum_1^n\log 1/\l_k}.
	$$
 \end{Prop}
 \begin{proof} 
     It is not restrictive to assume $a=0$, $b=1$.  Then the
     eigenvalues of $|D_{f}|^{-1}$ are given by $\La_{k} =
     \prod_{j=0}^{k}\l_{j}$, with multiplicity $2P_{k} = 2
     \prod_{j=0}^{k}p_{j}$, $k\in\bn\cup\{0\}$, where we have set
     $\l_{0}:=1$, $p_{0}:=1$.  Therefore,
	 \begin{align*}
	     Tr(|D_{f}|^{-\a}) & 
	     = 2\sum_{k=0}^{\infty} \prod_{i=1}^{k} p_{i}(\l_{i})^{\a} 
	     = 2\sum_{k=0}^{\infty} P_{k}\La_{k}^{\a} \\
	     & = 2\sum_{k=0}^{\infty} \exp\left(\log P_{k}
		\left(1-\a\frac{\log 1/\La_{k}}{\log P_{k}}\right)\right).
	 \end{align*}
	 As in the proof of the previous Theorem we conclude.
 \end{proof}

 \begin{rem}
	In \cite{GuIs11} we define pointwise tangential upper and
	lower dimensions for subspaces of $\br^{n}$.  It turns out
	that for the uniformly generated symmetric fractals, such
	dimensions are constant and equal respectively to the maximal
	and minimal dimension computed below.
 \end{rem}

 
 \begin{Thm}\label{Thm:delta*symm}
	Let $(\ca,\ch,D)$ be the spectral triple associated with a
	uniformly generated symmetric fractal $F$, where the
	similarities $w_{n,i}$, $i=1,\dots,p_{n}$ have scaling
	parameter $\l_{n}$.  Then
	\begin{align*}
		\subc(\ca,\ch,D)
		&=\liminf_{n,k}\frac{\sum_{j=n}^{n+k}\log p_{j}}
		{\sum_{j=n}^{n+k}\log 1/\l_j},\cr
		\supc(\ca,\ch,D)
		&=\limsup_{n,k}\frac{\sum_{j=n}^{n+k}\log p_{j}}
		{\sum_{j=n}^{n+k}\log 1/\l_j}.
	\end{align*}
 \end{Thm}

 As before, we first discuss the lacunary case.

 \begin{Prop}\label{Thm:delta*symmL}
	Let $F$ be a uniformly generated symmetric fractal, with the
	notations above.  Then
	\begin{align*}
	\subc(\ca,\ch_{\ell},D_{\ell})
	&=\liminf_{n,k}\frac{\sum_{j=n}^{n+k}\log p_{j}}{\sum_{j=n}^{n+k}\log 
	1/\l_j},\cr
	\supc(\ca,\ch_{\ell},D_{\ell})
	&=\limsup_{n,k}\frac{\sum_{j=n}^{n+k}\log p_{j}}{\sum_{j=n}^{n+k}\log 
	1/\l_j}.
	\end{align*}
 \end{Prop}

 \begin{proof}
     Making use of the definitions in (\ref{eigenvalues}), one gets
     $$
     \m_{|D_{\ell}|^{-1}}(x)=\tilde\La_{k},\quad\sum_{m=0}^{k-1}\tilde P_{m}
     < x \leq \sum_{m=0}^{k}\tilde P_{m}.
     $$
     Because of Lemma \ref{Lemma:doublelim}, $\subc^{-1}$, resp. 
     $\supc^{-1}$, is equal to the $\limsup$, resp.  $\liminf$ when
     $t$ and $h$ go to $\infty$, of the quantity $\frac{1}{h}(\log
     1/\m(e^{t+h})-\log1/\m(e^{t}))$, which may be rewritten as
     \begin{equation}\label{ratio}
	\frac{\log 1/\tilde\La_{k}-\log 1/\tilde\La_{m}}
	{\log\left(\sum_{j=0}^{k}\tilde P_{j}-\th_{k}\tilde P_{k}\right)
	-\log\left(\sum_{j=0}^{m}\tilde P_{j}-\th'_{m}\tilde P_{m}\right)}
    \end{equation}
    for suitable constants $\th_{k}$, $\th'_{k}$ in $[0,1)$.  Since
    the denominator goes to infinity, additive perturbations of the
    numerator and of the denominator by bounded sequences do not alter
    the $\limsup$, resp.  $\liminf$, therefore the uniform generation
    hypotheses imply that the ratio (\ref{ratio}) can be replaced by
    \begin{equation}\label{ratio2}
	\frac{\log 1/\La_{k}-\log 1/\La_{m}}
	{\log P_{k}-\log P_{m}}.
    \end{equation}
    Finally, since the denominator $\log P_{k}-\log P_{m}$ goes to
    infinity if and only if $k-m\to\infty$, the thesis follows.
 \end{proof}

 Then we discuss the filled case. Indeed for such spectral triple 
 the result holds in more generality, namely for any symmetric fractal.

 \begin{Prop}\label{Thm:delta*symmF}
	Let $F$ be a symmetric fractal, with the notations above.  Then
	\begin{align*}
	\subc(\ca,\ch_{f},D_{f})
	&=\liminf_{n,k}\frac{\sum_{j=n}^{n+k}\log p_{j}}{\sum_{j=n}^{n+k}\log 
	1/\l_j},\cr
	\supc(\ca,\ch_{f},D_{f})
	&=\limsup_{n,k}\frac{\sum_{j=n}^{n+k}\log p_{j}}{\sum_{j=n}^{n+k}\log 
	1/\l_j}.
	\end{align*}
 \end{Prop}

 \begin{proof}
	 In this case, the eigenvalues of $|D_{f}|^{-1}$ are the numbers
	 $\La_{k}$, each with multiplicity $2P_{k}$.  Therefore, the
	 quantity $\frac{1}{h}(\log 1/\m(e^{t+h})-\log1/\m(e^{t}))$,
	 may be rewritten as
     \begin{equation}\label{ratioF}
	\frac{\log 1/\La_{k}-\log 1/\La_{m}}
	{\log\left(\sum_{j=0}^{k} P_{j}-\th_{k} P_{k}\right)
	-\log\left(\sum_{j=0}^{m} P_{j}-\th'_{m} P_{m}\right)}
    \end{equation}
    for suitable constants $\th_{k}$, $\th'_{k}$ in $[0,1)$.  
	
	Let us observe that, since $p_{i}\geq2$, 
	\begin{align*}
		\log\left(\sum_{j=0}^{k} P_{j}-\th_{k} P_{k}\right) - \log P_{k}
		&\leq\log\frac{\sum_{j=0}^{k} P_{j}}{P_{k}}\\
		&=\log\left(\sum_{j=0}^{k}\prod_{i=j+1}^{k}\frac1{p_{i}} \right)
		\leq\log2,
	\end{align*}
	therefore, as before, the ratio above may be replaced by
    \begin{equation}\label{ratio2F}
	\frac{\log 1/\La_{k}-\log 1/\La_{m}}
	{\log P_{k}-\log P_{m}}.
    \end{equation}
    The thesis follows as before.
 \end{proof}
 
 Now we turn to the triple $(\ca,\ch,D)$.
 
 \begin{proof} {\it (of Theorem \ref{Thm:delta*symm}).}
	 The result will follow from Propositions \ref{Thm:delta*symmL}
	 and \ref{Thm:delta*symmF} if we show that the assumptions of Proposition 
	 \ref{dim.tang.eq} are satisfied. Let us observe that
	 $$
	 A:= \max\left(
	 \sup_{k}\frac{\tilde\La_{k}}{\La_{k}},
	 \sup_{k}\frac{\La_{k}}{\tilde\La_{k}}\right)
	  =
	 \sup_{k}\frac{p_{k+1}-1}{1-p_{k+1}\l_{k+1}}
	 $$
	 is finite by hypothesis. Also
	 $$
	 B:= \max\left(
	 \sup_{m}\frac{\sum_{m=0}^{k+1}\tilde P_{m}}{\sum_{m=0}^{k} P_{m}},
	 \sup_{m}\frac{\sum_{m=0}^{k+1} P_{m}}{\sum_{m=0}^{k}\tilde P_{m}}
	 \right)
	 $$
	 is finite. Indeed, setting $p=\sup_{m}p_{m}$,
	 \begin{align*}
		 \frac{\sum_{m=0}^{k+1}\tilde P_{m}}{\sum_{m=0}^{k} P_{m}}
		 &=\frac{\sum_{m=0}^{k+1}(p_{m+1}-1) P_{m}}{\sum_{m=0}^{k} 
		P_{m}}\\
		 &\leq (p-1)\left(1+\frac{P_{k+1}}{\sum_{m=0}^{k} 
		P_{m}}\right)\\
		 &=(p-1)\left(1+\frac{1}{\sum_{m=0}^{k} \prod_{j=m+1}^{k+1}
		p_{j}^{-1}}\right)
		 \leq p^2-1,
	 \end{align*}
	 and the other bound is obtained in the same way.
	 Therefore, if $\sum_{m=0}^{k-1}\tilde P_{m}
     < x \leq \sum_{m=0}^{k}\tilde P_{m}$ and $\sum_{m=0}^{k-1} P_{m}
     < y \leq \sum_{m=0}^{k} P_{m}$, then $x/y\leq B$, hence
	 $$
	 \mu_{|D_{\ell}|^{-1}}(x)=\tilde{\Lambda}_{k}
	 \leq A \Lambda_{k} = A \mu_{|D_{f}|^{-1}}(y)
	 \leq  A \mu_{|D_{f}|^{-1}}\left(\frac{x}{B}\right).
	 $$
	 The inequality in the other direction is proved in the same way.
 \end{proof}

 
 We conclude this section with a corollary of the theorems above and 
 of Theorem \ref{h-Minkowski}
 
 \begin{Cor}
     Assume $F$ is a uniformly generated symmetric fractal, which is
     $h$-Minkowski measurable, $h\in G_{d}$.  Then $\subc = d = \supc
     = \supd_{B}(F)$.
 \end{Cor}
 
 \section{Appendix. H\"older inequalities for singular 
 traces}\label{Appendix}
 
 For the reader's convenience, we recall some notions from \cite{GuIs1}
 that will be needed in this section.  
 
 Let $\t_{\o}$ be a singular trace on $\cb(\ch)$.  Then there is a unique 
 positive linear functional $\f$ on the cone of positive 
 non-increasing right-continuous functions on $[0,\infty)$, which is 
 dilation-invariant ($i.e.$ $\a\f(D_{\a}\m) = \f(\m)$, where 
 $D_{\a}\m(t) := \m(\a t)$, $\a,t>0$) and such that $\t_{\o}(a) = 
 \f(\m_{a})$, for any positive compact operator $a$.  In particular 
 the domain of the singular trace consists of the elements $a$ for 
 which $\f(\m_{a})$ is finite.
 
 It is not known if every positive linear dilation-invariant 
 functional $\f$ gives rise to a singular trace $\t_{\o}$ on $\cb(\ch)$ via 
 the formula $\t_{\o}(a):= \f(\m_{a})$.  But this is true if the functional 
 $\f$ is monotone, $i.e.$ increasing (which means $\sinc_{\m} \leq 
 \sinc_{\n}$ implies $\f(\m)\leq \f(\n)$) or decreasing (which means 
 $\sdec_{\m}\leq\sdec_{\n}$ implies $\f(\m)\leq \f(\n)$).  This means 
 in particular that all known formulas for singular traces are given 
 by a monotone functional.
 
 We can now state the main result of this appendix.
 
 \begin{Thm}
 	For any $a,b$ in the domain of $\t_{\o}$, $p,q\in[1,\infty]$ conjugate 
 	exponents, there holds
	$$
	|\t_{\o}(ab)| \leq \t_{\o}(|ab|) \leq C_{p} \t_{\o}(|a|^{p})^{1/p}\ \t_{\o}(|b|^{q})^{1/q},
	$$
	where $C_{p} := 1+ 2\frac{\sqrt{p-1}}{p}\in[1,2]$. \\
	If $\t_{\o}$ is generated by a monotone $\f$, then one can choose 
	$C_{p}=1$, for any $p\in[1,\infty]$.
 \end{Thm}
 \begin{proof}
     It is a consequence of the following propositions.
 \end{proof}
 
 \begin{Prop}
     For any $a,b$ in the domain of $\t_{\o}$, $p,q\in[1,\infty]$ conjugate
     exponents, there holds
     $$
	|\t_{\o}(ab)| \leq \t_{\o}(|ab|) \leq C_{p} \t_{\o}(|a|^{p})^{1/p}\ \t_{\o}(|b|^{q})^{1/q},
     $$
     where $C_{p} := 1+ 2\frac{\sqrt{p-1}}{p}\in[1,2]$.
 \end{Prop}	
 \begin{proof}
	Let $\a,\b>0$. Then, for any $t>0$, one has $\m_{|ab|}((\a+\b)t) 
	\leq \m_{a}(\a t)\m_{b}(\b t)$, $i.e.$, $D_{\a+\b}\m_{|ab|} \leq 
	D_{\a}\m_{a}D_{\b}\m_{b}$, so that
	\begin{align*}
		\t_{\o}(|ab|) & = \f(\m_{|ab|}) = (\a+\b)\f(D_{\a+\b}\m_{|ab|}) \\
		& \leq (\a+\b) \f(D_{\a}\m_{a}D_{\b}\m_{b}) \\
		& \leq (\a+\b) \left\{ \frac{1}{p} 
		\f((D_{\a}\m_{a})^{p}) + \frac{1}{q} 
		\f((D_{\b}\m_{b})^{q}) \right\}\\
		& = (\a+\b) \left\{ \frac{1}{\a p} 
		\f(\m_{a}^{p}) + \frac{1}{\b q} 
		\f(\m_{b}^{q}) \right\} \\
		& = (\a+\b) \left\{ \frac{1}{\a p} 
		\t_{\o}(|a|^{p}) + \frac{1}{\b q} 
		\t_{\o}(|b|^{q}) \right\},
	\end{align*}
	where we used Young's inequality and the dilation invariance of $\f$.
	Therefore, substituting $a/\t_{\o}(|a|^{p})^{1/p}$ for $a$, and 
	$b/\t_{\o}(|b|^{q})^{1/q}$ for $b$, we get
	$$
	\t_{\o}(|ab|) \leq \frac{(\a+\b)(\a p + \b q)}{\a\b pq} 
	\t_{\o}(|a|^{p})^{1/p} \t_{\o}(|b|^{q})^{1/q}.
	$$
	Set $g(x) := 1+ \frac{x}{p} + \frac{1}{xq}$, so that $g(\b/\a) = 
	\frac{(\a+\b)(\a p + \b q)}{\a\b pq}$. Minimizing $g$ over 
	$(0,\infty)$, we obtain $\min_{\a,\b>0} \frac{(\a+\b)(\a p + \b 
	q)}{\a\b pq} = 1+ 2\frac{\sqrt{p-1}}{p}$, which is easily seen to 
	belong to $[1,2]$.	
 \end{proof}
 
 For monotone $\f$'s the result is contained in the following
 propositions, but we need a preliminary result, which is interesting on
 its own.
 
 \begin{Prop}\label{prop:coWeyl}
 	Let $a,b\ \in\ck(\ch)$. Then, for any $n\in\bn\cup\{0\}$, 
	$$
	\sum_{i=2n}^{\infty} \m_{i}(ab) \leq \sum_{i=n}^{\infty} 
	\m_{i}(a)\m_{i}(b).
	$$
 \end{Prop}
 \begin{proof}
 	Let us first assume that $a,b\geq 0$, and let $ab = v|ab|$ be the 
 	polar decomposition of $ab$. Let $\ch'$ be a separable infinite 
 	dimensional Hilbert space, and let $v'\in \cb(\ch\oplus\ch')$ be a 
 	partial isometry with initial space $\ker |ab| \oplus \ch'$ and 
 	final space $(\ran |ab|)^{\perp} \oplus \ch'$. Set $u^{*} := 
 	(v\oplus 0) + v'$, which is a unitary operator of $\ch\oplus\ch'$. 
 	Then $(a\oplus 0)(b\oplus 0) = u^{*}|(a\oplus 0)(b\oplus 0)|$. As 
 	$\m_{n}(a\oplus 0) = \m_{n}(a)$, as long as they are nonzero, we 
 	can assume that $|ab| = uab$, with a unitary operator $u$.
	
	Let us now recall that 
	\begin{align*}
		\sdec_{2n}(ab) & := \sum_{i=2n}^{\infty} \m_{i}(ab) = 
		\inf_{p\in\cp_{2n}} \Tr(|ab|p^{\perp}) \\
	    & = \inf_{p,q\in\cp_{n}} \Tr((p\vee q)^{\perp}uab(p\vee 
	    q)^{\perp}),
	\end{align*}
	where $\cp_{n} := \{ p\in Proj(\ch\oplus\ch') : p\leq 
	\ker(ab)^{\perp},\ \Tr(p) \leq n\}$. Let us denote by $e_{n}(x)$ the 
	orthogonal projection on the eigenspace of the compact operator $x$ 
	associated to the largest $n$ eigenvalues. Choose $p := 
	e_{n}(uau^{*})\in\cp_{n}$ and $q := e_{n}(b)\in\cp_{n}$, and denote 
	by $a_{n}:= ae_{n}(a)^{\perp}$, $b_{n} := bq^{\perp}$, so that 
	$p^{\perp}uau^{*} = ua_{n}u^{*}$. Then
	\begin{align*}
		\sdec_{2n}(ab) & \leq \Tr((p\vee q)^{\perp} p^{\perp} uau^{*}ub 
		q^{\perp}(p\vee q)^{\perp}) \\
		& = \Tr((p\vee q)^{\perp} ua_{n}u^{*}ub_{n}(p\vee q)^{\perp}) \\
		& = |\Tr((p\vee q)^{\perp} ua_{n}b_{n})| 
		 \leq \Tr(|(p\vee q)^{\perp} ua_{n}b_{n}|) \\
		& \leq \Tr(|ua_{n}b_{n}|) = \Tr(|a_{n}b_{n}|) \\
		& = \sum_{i=0}^{\infty} \m_{i}(a_{n}b_{n}) 
		 \leq \sum_{i=0}^{\infty} \m_{i}(a_{n}) \m_{i}(b_{n}) \\
		& = \sum_{i=0}^{\infty} \m_{n+i}(a) \m_{n+i}(b),
	\end{align*}
	where the last inequality is Weyl's inequality (\cite{GK}, page 49).
	
	In the case where $a,b$ are arbitrary compact operators, let 
	$a = u|a|$, $b = |b^{*}|v$, $ab = w|ab|$ be polar decompositions. Then
	\begin{align*}
		\sdec_{2n}(ab) & = \sdec_{2n}(|ab|) = \sdec_{2n}(w^{*}ab) = 
		\sdec_{2n}(w^{*}u|a||b^{*}|v) \\
		& \leq \|w^{*}u\| \|v\| \sdec_{2n}(|a||b^{*}|) \\
		& \leq \sdec_{2n}(|a||b^{*}|) \\
		& \leq \sum_{i=n}^{\infty} \m_{i}(|a|) \m_{i}(|b^{*}|) \\
		& = \sum_{i=n}^{\infty} \m_{i}(a) \m_{i}(b),
	\end{align*}
	where we used $\m_{i}(|x|)=\m_{i}(x)=\m_{i}(x^{*})$, and, in the last 
	inequality, the thesis already proved for positive compact operators.
 \end{proof}
 
 \begin{Prop}
	Let $\f$ be a positive linear dilation invariant monotone 
	functional, and let $\t_{\o}(a) = \f(\m_{a})$ be the singular trace it 
	defines.  Then, for any $a,b$ in the domain of $\t_{\o}$, 
	$p,q\in[1,\infty]$ conjugate exponents, there holds
	$$
	|\t_{\o}(ab)| \leq \t_{\o}(|ab|) \leq \t_{\o}(|a|^{p})^{1/p}\ \t_{\o}(|b|^{q})^{1/q}.
	$$
 \end{Prop}
 \begin{proof}
     $(i)$ Let us first assume that $\f$ is increasing.
	Let us introduce the functions $\m_{a}(t) := \m_{n}(a), t\in[n,n+1)$, 
	and analogously for $\m_{b}$. Then Weyl's inequality reads as follows
	$$
	\int_{0}^{t} \m_{ab}(s)ds \leq \int_{0}^{t} 
	\m_{a}(s)\m_{b}(s) ds.
	$$
	As $\f$ is increasing, we get
	\begin{align*}
		\t_{\o}(|ab|) & = \f(\m_{ab}) \leq \f(\m_{a}\m_{b}) \\
		& \leq \frac{1}{p} \f(\m_{a^{p}}) + \frac{1}{q} \f(\m_{b^{q}}) \\
		& = \frac{1}{p} \t_{\o}(a^{p}) + \frac{1}{q} \t_{\o}(b^{q}),
	\end{align*}
	where we used Young's inequality for real numbers and the properties 
	of $\f$.
	Therefore, substituting $a/\t_{\o}(|a|^{p})^{1/p}$ for $a$, and 
	$b/\t_{\o}(|b|^{q})^{1/q}$ for $b$, we get
	$$
	\t_{\o}(|ab|) \leq \t_{\o}(|a|^{p})^{1/p} \t_{\o}(|b|^{q})^{1/q}.
	$$
 $(ii)$ Assume now that $\f$ is decreasing.
	Then the inequality in 
	Proposition \ref{prop:coWeyl} reads as follows
	$$
	2\int_{t}^{\infty} \m_{ab}(2s)ds = \int_{2t}^{\infty} \m_{ab}(s)ds 
	\leq \int_{t}^{\infty} \m_{a}(s)\m_{b}(s) ds.
	$$
	As $\f$ is decreasing and dilation invariant, we get
	\begin{align*}
		\t_{\o}(|ab|) & = \f(\m_{ab}) = 2\f(D_{2}\m_{ab}) \leq \f(\m_{a}\m_{b}) \\
		& \leq \frac{1}{p} \f(\m_{a^{p}}) + \frac{1}{q} \f(\m_{b^{q}}) \\
		& = \frac{1}{p} \t_{\o}(a^{p}) + \frac{1}{q} \t_{\o}(b^{q}),
	\end{align*}
	where we used Young's inequality for real numbers and the properties 
	of $\f$.
	Therefore, the thesis follows as in $(i)$. 	
 \end{proof}

\begin{ack} We thank Michel L. Lapidus for discussions, and for 
pointing out the relevance of some results in \cite{HeLa} for our 
lines of research. We also thank Alain Connes for suggestions.
\end{ack} 



\begin{thebibliography}{99}

 \bibitem{AGPS} S. Albeverio, D.~Guido, A.~Ponosov, S.~Scarlatti.
 {\it Singular traces and compact operators}. J. Funct. Anal., 
 {\bf 137} (1996), 281--302.
 
 \bibitem{BeFa} M-T. Benameur, T. Fack.  {\it On von Neumann spectral 
 triples}, preprint math.KT/0012233

 \bibitem{BS} C. Bennett, R. Sharpley. {\it Interpolation of
 operators}. Academic Press, New York, 1988.
 
 \bibitem{BGT} N.H. Bingham, C.M. Goldie , J.L. Teugels {\it Regular 
 variation} Cambridge University Press, Cambridge, 1987.
 
 \bibitem{CPS} A.L. Carey, J. Phillips, F. Sukochev
 {\it Spectral flow and Dixmier traces}, preprint math.OA/0205076 

 \bibitem{Co} A. Connes. {\it Non Commutative Geometry}. Academic
 Press, 1994.
 
 \bibitem{Cpc} A. Connes, private communication.
 
 \bibitem{CiGS1} F. Cipriani, D. Guido, S. Scarlatti.  \textit{ A
 remark on trace properties of K-cycles}.  J. Operator Theory,
 \textbf{ 35} (1996), 179--189.

 \bibitem{Dixmier} J. Dixmier.  {\it Existence de traces non
 normales}.  C.R. Acad.  Sci.  Paris, {\bf 262} (1966), 1107--1108.
 
  \bibitem{DPSS} P.G. Dodds, B. de Pagter, E.M. Semenov, F. A. 
  Sukochev.  {\it Symmetric functionals and singular traces}.  
  Positivity {\bf 2} (1998),  47--75.
 
 \bibitem{DFWW} K. Dykema, T. Figiel, G. Weiss, M. Wodzicki.  {\it 
 Commutator Structure of Operator Ideals}, MSRI Preprint n. 2001-013

 \bibitem{FK} T. Fack, H. Kosaki.  {\it Generalized s-numbers of
 $\t$-measurable operators}.  Pacific J. Math., {\bf 123} (1986), 269.

 \bibitem{Fur} H. Furstenberg. Talk at the Conference ``Fractals in 
 Graz'', June 2001.
 
 \bibitem{Gatzuras} D. Gatzouras.  \textit{ Lacunarity of self-similar
 and stochastically self-similar sets} Trans.  Amer.  Math.  Soc.,
 \textbf{ 352} (2000), 1953--1983.

 \bibitem{GK} I. Gohberg, M.G. Krein. {\it Introduction to the
 theory of non-selfadjoint operators}, Mos\-cow (1985).

 \bibitem{GBVa} J.M. Gracia-Bondia, J.C. V\'arilly, H. Figueroa.  {\it 
 Elements of noncommutative geometry}.  Birkh\"auser Boston, Inc., 
 Boston, MA, 2001.
 
 \bibitem{GuIs1} D. Guido, T. Isola.  {\it Singular traces for 
 semifinite von~Neumann algebras}.  Journal of Functional Analysis, 
 {\bf 134} (1995), 451--485.

 \bibitem{GuIs4} D. Guido, T. Isola.  {\it Noncommutative Riemann 
 integration and Novikov-Shubin invariants for open manifolds}.  
 J. Funct. Anal. {\bf 176} (2000) 115--152

 \bibitem{GuIs5} D. Guido, T. Isola.  {\it Singular traces and
 Novikov-Shubin invariants}, in Operator Theoretical Methods,
 Proceedings of the 17th Conference on Operator Theory Timisoara
 (Romania) June 23-26, 1998, A. Gheondea, R.N. Gologan, D. Timotin
 Editors, The Theta Foundation, Bucharest 2000.

 \bibitem{GuIs8} D. Guido, T. Isola.  {\it Fractals in Noncommutative 
 Geometry}, in the Proceedings of the Conference "Mathematical Physics 
 in Mathematics and Physics", Siena 2000, Edited by Roberto Longo, 
 University of Rome II, Italy, Fields Institute Communications, Vol.  
 30 American Mathematical Society, Providence, RI, 2001


 \bibitem{GuIs10} D. Guido, T. Isola.  {\it Dimensions, spectral
 triples and fractals: fractals in $\br^{n}$}.  In progress.

 \bibitem{GuIs11} D. Guido, T. Isola.  {\it Tangential dimensions for 
 metric spaces and measures}. In progress.
 
 \bibitem{GuIs12} D. Guido, T. Isola.  D.~Guido, T.~Isola.  On the
 domain of singular traces.  {\it Int.  J. Math.}, {\bf 13} (2002),
 667-674.

 \bibitem{Taub} D. Guido, T. Isola.  In progress.
 
 \bibitem{HeLa} C.Q. He, M.L. Lapidus.  {\it Generalized Minkowski
 content, spectrum of fractal drums, fractal strings and the Riemann
 zeta-function}.  Mem.  Amer.  Math.  Soc.,  {\bf 127} (1997), no.  608.

 \bibitem{Hutch} J.E. Hutchinson. {\it Fractals and 
 self-similarity},  Indiana Univ.  Math.  J. \textbf{30} (1981),  
 713--747.
 
 \bibitem{KaSa} J-P. Kahane, R. Salem. \textit{ Ensembles parfaits et s\'eries 
 trigonom\'etriques}. Hermann, Paris, 1994.
 
 \bibitem{KiLa2} J. Kigami, M.L. Lapidus.  {\it Self-similarity of 
 volume measures for Laplacians on p.c.f. self-similar fractals}.  
 Comm.  Math.  Phys.  {\bf 217} (2001), 165--180.
 
\bibitem{LavF} M.L. Lapidus, M. van Frankenhuysen.  {\it Fractal 
geometry and number theory.  Complex dimensions of fractal strings and 
zeros of zeta functions}.  Birkh\"auser Boston, Inc., Boston, MA, 2000.

\bibitem{LaPo} M.L. Lapidus, C. Pomerance.  \textit{ The Riemann 
zeta-function and the one-dimensional Weyl-Berry conjecture for 
fractal drums}.  Proc.  London Math.  Soc.  (3) \textbf{ 66} (1993), 
41--69.

\bibitem{random} R.D. Mauldin, S.C. Williams.  {\it Random recursive 
constructions: asymptotic geometric and topological properties}.  
Trans.  Amer.  Math.  Soc.  {\bf 295} (1986), 325--346.
 

 \bibitem{Tricot} C. Tricot.  \textit{Curves and fractal dimension}.  
 Springer-Verlag, New York, 1995.

 \bibitem{Varga} J.V. Varga. {\it Traces on irregular ideals}. Proc.
 A.M.S., {\bf107} (1989), 715.

 

\end{thebibliography}
\end{document}